# TREE BASED FUNCTIONAL EXPANSIONS FOR FEYNMAN–KAC PARTICLE MODELS

By Pierre Del Moral, Frédéric Patras[1] and Sylvain Rubenthaler

*INRIA Centre Bordeaux Sud Quest, Université de Nice-Sophia Antipolis and Université de Nice-Sophia Antipolis*

We design exact polynomial expansions of a class of Feynman–Kac particle distributions. These expansions are finite and are parametrized by coalescent trees and other related combinatorial quantities. The accuracy of the expansions at any order is related naturally to the number of coalescences of the trees. Our results include an extension of the Wick product formula to interacting particle systems. They also provide refined nonasymptotic propagation of chaos-type properties, as well as sharp $\mathbb{L}_p$-mean error bounds, and laws of large numbers for $U$-statistics.

**0. Introduction.** The typical phenomenon we are interested in is the behavior of particles scattered according to a diffusion process (encoded by a Markov transition) and submitted to a potential function. The problem is fairly general, since its modelization encodes, besides the scattering of particles in an absorbing medium, also filtering problems in signal processing or polymerization models, to quote only a few of the applications areas of the theory. In these three cases, the potential function represents, respectively, the absorption rate, the likehood of the signal values conditional to some observations, or the intermolecular attraction or repulsion forces between the monomers. For further details and a survey of the applications of Feynman–Kac interacting particle models, the interested reader is recommended to consult the pair of books [5, 6], and the references therein.

As it is well known, naive Monte-Carlo methods fail to give a satisfactory answer to these questions. This is easily understood with the example of

Received November 2006; revised October 2008.
[1]Supported by ANR Grant AHBE 05-42234.

particles in an absorbing medium. In that case, the size of the nonabsorbed population will decrease according to the absorption rate of the medium, so that, in the end, the empirical distribution of the surviving sample may give only a poor approximation of the final repartition of the particles. The solution to the decrease of the population size is simple and relies on a mean-field approximation. In general, particles will explore the state space as a free Markov evolution; during their exploration particles with low potential are killed, while the ones with high potential value duplicate.

During the last two decades, the asymptotic analysis of these Feynman–Kac interacting particle models has been developed in various directions, including propagation of chaos analysis, $\mathbb{L}_p$-mean error estimates, central limit-type theorems, and large deviation principles. The purpose of the present work is to improve on these results and develop exact, nonasymptotic, tree-based functional representations of particle block distributions. Actually, two probability distributions arise naturally in these models. To stick with the diffusion-absorption scheme, they are associated, respectively, to the empirical distribution of the nonabsorbed component of the Markov chain modeling the trajectories (1.8) and to the corresponding normalized empirical occupation measure (1.9). In both cases, the distributions can be expanded polynomially with respect to $N^{-1}$, where $N$ stands for the total population size. These Laurent-type expansions are the main results of the article. They are described in Theorems 1.5 and 4.12.

They rely on an original combinatorial, and permutation group analysis on a special class of trees and forests that parametrize naturally the trajectories of interacting particle systems. Much attention has been paid recently to the combinatorics of trees and their applications in physics. A short overview of the areas where these phenomena show up is included in Section 3. To the best of our knowledge, their introduction in the analysis of Feynman–Kac and interacting particle models is new. Since the corresponding geometrical and combinatorial study might be useful in other contexts, we have tried to isolate its exposition, as far as possible, from its applications to Feynman–Kac models.

The functional expansions allow us to derive as a direct consequence Wick formulas for interacting particle systems [Theorem 3.14, (5.1)], refined, nonasymptotic, propagation of chaos-type properties including strong expansions of the particle block distributions with respect to Zolotarev-type seminorms (5.2), as well as explicit formulas for $\mathbb{L}_p$-mean error bounds (5.4), and laws of large numbers for $U$-statistics (5.3), yielding what seems to be the first results of this type for mean-field particle models.

The article is divided into four main parts, devoted, respectively, to the precise description of Feynman–Kac particle models and tree-based expansions (Section 1), to the proof of the expansion formulas (Section 2), to the combinatorial and group-theoretical analysis of trajectories in mean-field interacting particle models (Section 3), and to the extension to path-space models with applications to propagation of chaos-type properties (Section 4). The last section gathers various direct consequences of the main theorems.

*Notation.*

*Combinatorial quantities.* Let us gather first various notations, for further use in the various parts of the article.

Let $S$ be a finite set; we write $|S|$ for the cardinal of $S$. The notation is extended to applications $a$ between finite sets and $|a|$ denotes then the cardinal of the image of $a$. Precisely, for any pair of integers $(l,m) \in (\mathbb{N}^\star)^2$, we set $[l] = \{1, \ldots, l\}$, and $[l]^{[m]}$ the set of mappings $a$ from $[m]$ into $[l]$. By $|a|$, we denote the cardinality of the set $a([m])$, and for any $1 \leq p \leq l$ we set

$$[l]^{[m]}_p = \{a \in [l]^{[m]} : |a| = p\}.$$

When $l > m$, we also denote by $\langle m, l \rangle (= [l]^{[m]}_m)$ the set of all

$$(l)_m := \frac{l!}{(l-m)!}$$

one-to-one mappings from $[m]$ into $[l]$, and by $\mathbf{S}_l = \langle l, l \rangle$ the symmetric group of all permutations of $[l]$. We also denote by $1_l$ the identity in $\mathbf{S}_l$. Recall that the Stirling number of the second kind $S(q,p)$ is the number of partitions of $[q]$ into $p$ nonempty subsets, so that

$$|[N]^{[q]}_p| = S(q,p)(N)_p \quad \text{and} \quad N^q = \sum_{1 \leq k \leq q} S(q,k)(N)_k.$$

Sequences of maps and associated quantities will be one of the main ingredients of our approach to Feynman–Kac models. We write $\mathcal{A}_{n,q} =_{\text{def.}} ([q]^{[q]})^{n+1}$ (resp. $\mathcal{IA}_{n,q}$) for the set of $(n+1)$-sequences of maps (resp. of weakly increasing maps) $\mathbf{a} = (a_p)_{0 \leq p \leq n}$ from $[q]$ into itself:

$$[q] \xleftarrow{a_0} [q] \xleftarrow{a_1} \cdots \longleftarrow [q] \xleftarrow{a_{n-1}} [q] \xleftarrow{a_n} [q].$$

Notice that we write in bold the symbols for sequences (of maps, integers, ...) such as $\mathbf{a}$.

For any sequences of integers $\mathbf{p} = (p_k)_{0 \leq k \leq n}$, and $\mathbf{l} = (l_k)_{0 \leq k \leq n}$, we write $\mathbf{p} \leq \mathbf{l}$ if and only if $p_k \leq l_k$ for all $0 \leq k \leq n$. We write $\|\mathbf{p}\|$ for $(p_0 + \cdots + p_n)$. Assuming now that $\mathbf{p} \leq \mathbf{l}$, we use the multi-index notation

$$(\mathbf{l})_\mathbf{p} = \prod_{0 \leq k \leq n} (l_k)_{p_k}, \qquad \mathbf{p}! = \prod_{0 \leq k \leq n} p_k!, \qquad s(\mathbf{l}, \mathbf{p}) = \prod_{k=0}^n s(l_k, p_k),$$

where the $s(l_k, p_k)$ are Stirling numbers of the first kind. Recall that these numbers provide the coefficients of the polynomial expansion of the $(N)_p$ (Stirling formula):

$$(N)_p = \sum_{1 \leq k \leq p} s(p,k) N^k;$$

see, for example, [3] for further details. The difference $(\mathbf{p} - \mathbf{l})$ and, respectively, the addition $(\mathbf{p} + \mathbf{l})$ of two sequences is the sequence $(p_k - l_k)_{0 \leq k \leq n}$ and, respectively, $(p_k + l_k)_{0 \leq k \leq n}$. When no confusions can arise, we write $\mathbf{N}$ and $\mathbf{q}$, for the constant sequences $(N)_{0 \leq i \leq n}$, and $(q)_{0 \leq i \leq n}$. We also write $\mathbf{1}$ and, respectively, $\mathbf{0}$, for the sequence of unit integers and, respectively, null integers. The above definitions are extended to infinite sequence of integers $\mathbf{p} = (p_k)_{k \geq 0} \in \mathbb{N}^{\mathbb{N}}$, with a finite number of strictly positive terms. Any function $\alpha : \mathbb{N} \mapsto \mathbb{N}$ on the set of integers into itself, is extended to integer sequences $\mathbf{p} \in \mathbb{N}^{\mathbb{N}}$, by setting $\alpha(\mathbf{p}) = (\alpha(p_k))_{k \geq 0}$.

In particular, for $\mathbf{a} \in \mathcal{A}_{n,q}$, we write $\|\mathbf{a}\|$ the sequence $(\|a_i\|)_{0 \leq i \leq n}$. By a slight abuse of notation, we write $\|\mathbf{a}\|$ for $\|(|\mathbf{a}|)\| = |a_0| + \cdots + |a_n|$. The sequence $\mathbf{q} - |\mathbf{a}|$ is called the coalescence sequence of $\mathbf{a}$. The coalescence degree of $\mathbf{a}$ is defined by

$$coal(\mathbf{a}) := \|\mathbf{q}\| - \|\mathbf{a}\| = \sum_{i=0}^{n} (q - |a_i|).$$

The subset of $\mathcal{A}_{n,q}$ of sequences $\mathbf{a}$ such that $|\mathbf{a}| \geq \mathbf{q} - \mathbf{p}$ is written $\mathcal{A}_{n,q}(\mathbf{p})$.

*Measures, norms and related notions.* Let $(E, \mathcal{E})$ be a measurable space. We denote, respectively, by $\mathcal{M}(E)$, $\mathcal{P}(E)$ and $\mathcal{B}_b(E)$, the set of all finite signed measures on $(E, \mathcal{E})$, the convex subset of all probability measures and the Banach space of all bounded and measurable functions $f$ on $E$, equipped with the uniform norm $\|f\| = \sup_{x \in E} |f(x)|$. The total variation norm is written $\|\cdot\|_{\mathrm{TV}}$, so that, for any linear operator $L$ on $\mathcal{B}_b(E)$,

$$\|L\|_{\mathrm{TV}} := \sup_{f \in \mathcal{B}_b(E) : \|f\| \leq 1} |L(f)|.$$

A bounded integral operator $Q$ between the measurable spaces $(E, \mathcal{E})$ and $(F, \mathcal{F})$, such that, for any $f \in \mathcal{B}_b(F)$, the functions

$$Q(f) : x \in E \mapsto Q(f)(x) = \int_F Q(x, dy) f(y) \in \mathbb{R}$$

are $\mathcal{E}$-measurable, and bounded, generates a dual operator $\mu \mapsto \mu Q$ from $\mathcal{M}(E)$ into $\mathcal{M}(F)$, defined by $(\mu Q)(f) := \mu(Q(f))$. The tensor power $Q^{\otimes q}$ represents the bounded integral operator defined for any $f \in \mathcal{B}_b(F^q)$ by

$$Q^{\otimes q}(f)(x_1, \ldots, x_q) = \int_{F^q} [Q(x_1, dy_1) \cdots Q(x_q, dy_q)] f(y_1, \ldots, y_q).$$

For a bounded integral operator $Q_1$ from $E$ into $F$, and an operator $Q_2$ from $F$ into $G$, we denote by $Q_1Q_2$ the composition operator from $E$ into $G$, defined for any $f \in \mathcal{B}_b(G)$ by $(Q_1Q_2)(f) := Q_1(Q_2(f))$.

The notion of differential for sequences of signed measures is also useful. Let $(\Theta^N)_{N\geq 1} \in \mathcal{M}(E)^{\mathbb{N}}$ be a uniformly bounded sequence of signed measures on a measurable space $(E, \mathcal{E})$, in the sense that $\sup_{N\geq 1} \|\Theta^N\|_{\mathrm{TV}} < \infty$. The sequence $\Theta^N$ is said to converge strongly to some measure $\Theta \in \mathcal{M}(E)$, as $N \uparrow \infty$, if and only if

$$\forall f \in \mathcal{B}_b(E) \qquad \lim_{N\uparrow\infty} \Theta^N(f) = \Theta(f).$$

DEFINITION 0.1. Let us assume that $\Theta^N$ converges strongly to $\Theta$. The discrete derivative of the sequence $(\Theta^N)_{N\geq 1}$ is the sequence of measures $(\partial \Theta^N)_{N\geq 1}$ defined by

$$\partial \Theta^N := N[\Theta^N - \Theta].$$

We say that $\Theta^N$ is differentiable, if $\partial \Theta^N$ is uniformly bounded, and if it strongly converges to some measure $\partial \Theta \in \mathcal{M}(E)$, as $N \uparrow \infty$.

The discrete derivative $\partial \Theta^N$ of a differentiable sequence can itself be differentiable. In this situation, the derivative of the discrete derivative is called the second derivative and it is denoted by $\partial^2 \Theta = \partial(\partial\Theta)$, and so on.

A sequence $\Theta^N$ that is differentiable up to order $(k+1)$ has the following representation:

$$\Theta^N = \sum_{0\leq l\leq k} \frac{1}{N^l}\partial^l \Theta + \frac{1}{N^{k+1}}\partial^{k+1}\Theta^N$$

with $\sup_{N\geq 1} \|\partial^{k+1}\Theta^N\|_{\mathrm{TV}} < \infty$, and the convention $\partial^0 \Theta = \Theta$, for $l = 0$.

**1. Feynman–Kac semigroups.**

1.1. *Definitions.* We let $(E_n, \mathcal{E}_n)_{n\geq 0}$ be a collection of measurable state spaces. In some applications, the series $E_n$ will be constant or there will exist canonical isomorphisms between the various state spaces (think to the distribution of particles in a fixed absorbing medium), but it may also happen that the various state spaces are significantly different (think to path estimation and smoothing in signal processing or to polymerization sequences of monomers). For a survey of the various applications of the theory and examples of families of state spaces $E_n$ occurring in practice, we refer to the pair of books [5, 6], and the references therein.

We consider a distribution $\eta_0$ on $E_0$, a collection of Markov transitions $M_n(x_{n-1}, dx_n)$ from $E_{n-1}$ into $E_n$, and a collection of $\mathcal{E}_n$-measurable and bounded potential functions $G_n$ on the state spaces $E_n$. We will always assume that the potential functions are chosen such that

$$(1.1) \qquad 0 < \inf_{x_n \in E_n} G_n(x_n) \leq \sup_{x_n \in E_n} G_n(x_n) < \infty.$$

We associate to these objects the Feynman–Kac measures defined for any function $f_n \in \mathcal{B}_b(E_n)$ by the following formulas.

DEFINITION 1.1. The Feynman–Kac measure $\gamma_n$ and its normalization $\eta_n$ are defined by

$$(1.2) \quad \gamma_n(f_n) := \mathbb{E}\left[ f_n(X_n) \prod_{0 \leq k < n} G_k(X_k) \right] \quad \text{and} \quad \eta_n(f_n) := \gamma_n(f_n)/\gamma_n(1),$$

where $(X_n)_{n \geq 0}$ is a Markov chain, taking values in the state spaces $(E_n)_{n \geq 0}$, with initial distribution $\eta_0$ on $E_0$, and elementary transitions $M_n$ from $E_{n-1}$ into $E_n$.

By the Markov property and the multiplicative structure of (1.2), we check that the flow $(\eta_n)_{n \geq 0}$ satisfies the following equation:

$$(1.3) \qquad \eta_{n+1} = \Phi_{n+1}(\eta_n),$$

where the transformations $\Phi_{n+1} : \mathcal{P}(E_n) \to \mathcal{P}(E_{n+1})$ are defined for any pair $(\eta_n, f_{n+1}) \in (\mathcal{P}(E_n) \times \mathcal{B}_b(E_{n+1}))$ as follows:

$$\Phi_{n+1}(\eta_n)(f_{n+1}) = \frac{\eta_n(Q_{n+1}(f_{n+1}))}{\eta_n(Q_{n+1}(1))}$$

$$\text{with } Q_{n+1}(x_n, dx_{n+1}) = G_n(x_n) \times M_{n+1}(x_n, dx_{n+1}).$$

The measures $\gamma_n$ can be expressed in terms of the flow $(\eta_p)_{0 \leq p \leq n}$ with the formulas

$$(1.4) \qquad \gamma_n(f_n) = \eta_n(f_n) \times \gamma_n(1) \qquad \text{with } \gamma_n(1) = \prod_{0 \leq p < n} \eta_p(G_p)$$

for any $f_n \in \mathcal{B}_b(E_n)$.

To illustrate their meaning and motivate further the forthcoming developments, let us develop with some details an important application of these Feynman–Kac flows. Let us consider the distribution of a particle evolving in an absorbing medium, with obstacles related to potential functions $G_n$, taking values in $]0, 1]$. In this context, the particle $Y_n$ evolves according to two separate mechanisms. First, it moves from a site $y_{n-1} \in E_{n-1}$, to another $y_n \in E_n$ according to elementary transitions $M_n(y_{n-1}, dy_n)$. Then, it is absorbed with a probability $1 - G_n(y_n)$, and placed in an auxiliary cemetery state $Y_n = c$; otherwise it remains in the same site. If we let $T$ be the random absorption time, then it is not difficult to check that

$$(1.5) \quad \gamma_n(f_n) = \mathbb{E}[f_n(Y_n)1_{T \geq n}] \quad \text{and} \quad \eta_n(f_n) = \mathbb{E}[f_n(Y_n)|T \geq n].$$

In the above displayed formulas, we have used the convention that $f_n(c) = 0$, when $Y_n = c$. In time-homogeneous settings ($G_n = G$ and $M_n = M$), we have ([5], Proposition 12.4.1)

$$\gamma_n(1) = \mathbb{P}(T \geq n) \sim e^{-\lambda n}$$

as $n \to \infty$. The positive constant $\lambda$ is a measure of the strength and trapping effects of the obstacle. By a lemma of Varadhan's, $\lambda$ coincides with the logarithmic Lyapunov exponent $\lambda_0$ of the transition operator $Q(x, dy) = G(x)M(x, dy)$. Whenever it exists, the corresponding eigenfunction $h$ of $Q$ represents the ground state of the operator $Q$; see, for example, [5].

1.2. *Mean-field interacting particle models.* A natural mean-field particle model associated with the nonlinear Feyman–Kac flow (1.3) is the $E_n^N$-valued Markov chain $\xi_n^{(N)} = (\xi_n^{i,N})_{1 \leq i \leq N}$ with elementary transitions defined for any $F \in \mathcal{B}_b(E_n^N)$ by

$$\mathbb{E}(F(\xi_n^{(N)})|\xi_{n-1}^{(N)}) = \Phi_n(m(\xi_{n-1}^{(N)}))^{\otimes N}(F) \qquad \text{with } m(\xi_{n-1}^{(N)}) = \frac{1}{N} \sum_{i=1}^{N} \delta_{\xi_{n-1}^{i,N}}.$$

(1.6)

In other terms, given the configuration $\xi_{n-1}^{(N)}$ at rank $(n-1)$, the particle system $\xi_n^{(N)}$ at rank $n$ consists of $N$ independent and identically distributed random variables with common distribution $\Phi_n(m(\xi_{n-1}^{(N)}))$. The initial configuration $\xi_0^{(N)}$ consists of $N$ independent and identically distributed random variables with distribution $\eta_0$. Although the dependency of the model on $N$ is strong, due to the mean-field nature of the model, we abbreviate $(\xi_n^{i,N})_{1 \leq i \leq N}$ to $(\xi_n^i)_{1 \leq i \leq N}$ when $N$ is fixed, except when we want to emphasize explicitly the dependency on $N$ of the model, for example, as in the definition of the measure $\mathbb{P}_{n,q}^N$ below. The same observation will apply for other random variables, measures or functions in a self-explanatory way.

Notice that

$$\Phi_n(m(\xi_{n-1}^{(N)}))(dx_n) = \frac{1}{\sum_{j=1}^{N} G_{n-1}(\xi_{n-1}^j)} \left( \sum_{i=1}^{N} G_{n-1}(\xi_{n-1}^i) M_n(\xi_{n-1}^i, dx_n) \right)$$

so that the particle model evolves as a genetic-type model with proportional selections, and mutation transitions dictated by the pair of potential-transition $(G_{n-1}, M_n)$.

DEFINITION 1.2. The approximation measure $\gamma_n^N$ and the normalized approximation measure $\eta_n^N$ associated with the Feynman–Kac measures $\gamma_n$ and $\eta_n$ are defined by the empirical occupation measures

$$\eta_n^N := m(\xi_n^{(N)})$$

and for any $f_n \in \mathcal{B}_b(E_n)$,

(1.7) $\quad \gamma_n^N(f_n) := \eta_n^N(f_n) \times \gamma_n^N(1) \quad$ with $\gamma_n^N(1) := \prod_{0 \leq p < n} \eta_p^N(G_p).$

The distribution of particle blocks for any size $q \leq N$ is written $\mathbb{P}_{n,q}^N$:

(1.8) $\qquad \mathbb{P}_{n,q}^N := \text{Law}(\xi_n^{1,N}, \ldots, \xi_n^{q,N}) \in \mathcal{P}(E_n^q).$

For a rather complete asymptotic analysis of these measures we again refer the reader to [5], and references therein. In particular, for any $f_n \in \mathcal{B}_b(E_n)$, we have the following almost convergence results:

$$\lim_{N \to \infty} \eta_n^N(f_n) = \eta_n(f_n) \quad \text{and} \quad \lim_{N \to \infty} \gamma_n^N(f_n) = \gamma_n(f_n).$$

It is well known that $\gamma_n^N$ is an unbias approximation measure of $\gamma_n$, in the sense that

$$\mathbb{E}(\gamma_n^N(f_n)) = \gamma_n(f_n).$$

However, it turns out that for any $n \geq 1$

$$\mathbb{E}(\eta_n^N(f_n)) = \mathbb{E}(f_n(\xi_n^{1,N})) \neq \eta_n(f_n).$$

This means that the mean-field particle model is not an exact sampling algorithm of the distributions $\eta_n$. In practice, it is clearly important to analyze the bias of these quantities, and more generally the one of the distributions $\mathbb{P}_{n,q}^N$. From the pure mathematical point of view, the $\mathbb{P}_{n,q}^N$ are better understood when they are connected with the following measures.

DEFINITION 1.3. The tensor product occupation measures (resp. the restricted tensor product occupation measures) on $E_n^q$ are given by

$$(\eta_n^N)^{\otimes q} := \frac{1}{N^q} \sum_{a \in [N]^{[q]}} \delta_{(\xi_n^{a(1)}, \ldots, \xi_n^{a(q)})}$$

$$\left[ \text{resp. } (\eta_n^N)^{\odot q} := \frac{1}{(N)_q} \sum_{a \in \langle q, N \rangle} \delta_{(\xi_n^{a(1)}, \ldots, \xi_n^{a(q)})} \right].$$

The corresponding unnormalized tensor product occupation measures $(\gamma_n^N)^{\otimes q}$ [resp. restricted unnormalized tensor product occupation measures $(\gamma_n^N)^{\odot q}$] are defined for any $F \in \mathcal{B}_b(E_n^q)$ by

$$(\gamma_n^N)^{\otimes q}(F) := (\eta_n^N)^{\otimes q}(F) \times (\gamma_n^N(1))^q \quad \text{and}$$

$$(\gamma_n^N)^{\odot q}(F) := (\eta_n^N)^{\odot q}(F) \times (\gamma_n^N(1))^q.$$

The corresponding nonnegative measure $\mathbb{Q}_{n,q}^N$, indexed by the particle block sizes, on the product state spaces $E_n^q$, is defined for any $F \in \mathcal{B}_b(E_n^q)$ by

(1.9) $$\mathbb{Q}_{n,q}^N(F) := \mathbb{E}((\gamma_n^N)^{\otimes q}(F)).$$

Notice that these various tensor product occupation measures are symmetry-invariant by construction. That is, for any $F \in \mathcal{B}_b(E_n^q)$ and any $\sigma \in \mathbf{S}_q$,

$$(\eta_n^N)^{\otimes q}(F) = (\eta_n^N)^{\otimes q}(F \circ \sigma) \quad \text{and} \quad (\eta_n^N)^{\odot q}(F) = (\eta_n^N)^{\odot q}(F \circ \sigma),$$

where $\sigma$ acts by permutation on $E_n^q$. In particular, we may assume without restriction in our forthcoming computations on $q$-tensor product occupation measures, that $F$ is a symmetric function, that is,

$$F = F_{\text{sym}} := \frac{1}{q!} \sum_{\sigma \in \mathbf{S}_q} F \circ \sigma.$$

We write from now on, $\mathcal{B}_b^{\text{sym}}(E_n^q)$ for the set of all symmetric functions in $\mathcal{B}_b(E_n^q)$.

Notice also the following symmetry property, essential in view of all the forthcoming computations. Since, conditional to $\xi_{n-1}^{(N)}$, the $\xi_n^i$ are i.i.d., for any $F \in \mathcal{B}_b(E_n^q)$ and any $a, b \in \langle q, N \rangle^2$, we have

$$\mathbb{E}[F(\xi_n^{a(1)}, \ldots, \xi_n^{a(q)})] = \mathbb{E}[F(\xi_n^{b(1)}, \ldots, \xi_n^{b(q)})].$$

For instance, we have that

(1.10) $$\mathbb{E}((\eta_n^N)^{\odot q}(F)) = \mathbb{E}(F(\xi_n^1, \ldots, \xi_n^q)) = \mathbb{P}_{n,q}^N(F).$$

Let us return to the diffusion-absorption model introduced in the previous section. As a consequence of the results presented in [5], the sequence of distributions $\mathbb{Q}_{n,q}^N$ converges, as $N$ tends to infinity, to the distribution of $q$ nonabsorbed, and independent particles $(Y_n^i)_{1 \leq i \leq q}$ evolving in the original absorbing medium. That is, we have that

$$\lim_{N \to \infty} \mathbb{Q}_{n,q}^N(F) = \gamma_n^{\otimes q}(F) = \mathbb{E}(F(Y_n^1, \ldots, Y_n^q) 1_{T^1 \geq n} \cdots 1_{T^q \geq n})$$

for any $F \in \mathcal{B}_b(E_n^q)$, and where $T^i$ stands for the random absorption time sequence of the chain $(Y_k^i)_{k \geq 0}$, with $1 \leq i \leq q$.

This article is mainly concerned with explicit expansions of the deterministic measures $\mathbb{P}_{n,q}^N \in \mathcal{P}(E_n^q)$, and $\mathbb{Q}_{n,q}^N \in \mathcal{M}(E_n^q)$, with respect to the precision parameter $N$.

1.3. *Expansion formulas.* Recall that $\mathcal{A}_{n,q}$ stands for the set of $(n+1)$-sequences of mappings from $[q]$ into itself. We let $\Delta_{n,q}$ be the nonnegative measure-valued functional on $\mathcal{A}_{n,q}$ defined by

$$(1.11) \quad \Delta_{n,q} : \mathbf{a} \in \mathcal{A}_{n,q} \mapsto \Delta_{n,q}^{\mathbf{a}} := (\eta_0^{\otimes q} D_{a_0} Q_1^{\otimes q} D_{a_1} \cdots Q_n^{\otimes q} D_{a_n}) \in \mathcal{M}(E_n^q),$$

where the operator $D_b$ stands for the Markov transition from $E_n^q$ into itself, associated with a mapping $b \in [q]^{[q]}$, and defined by

$$D_b(F)(x_n^1, \ldots, x_n^q) := F(x_n^{b(1)}, \ldots, x_n^{b(q)})$$

for any $F \in \mathcal{B}_b(E_n^q)$, and $(x^1, \ldots, x^q) \in E_n^q$. Notice that

$$D_a D_b = D_{ab}$$

for any pair of mappings $(a,b) \in ([q]^{[q]})^2$. For $u$ a linear combination $\sum_{i \in I} \alpha_i a_i$ of elements of $[q]^{[q]}$, we extend the definition of $D$ by linearity and write

$$D_u := \sum_{i \in I} \alpha_i D_{a_i}.$$

Notice that the measures $\Delta_{n,q}^{\mathbf{a}}$ inherit a remarkable invariance property from their set-theoretic definition. Namely, let us introduce the natural left action of the group $\mathbf{S}_q^{n+2}$ on $\mathcal{A}_{n,q}$ defined for all $\mathbf{a} \in \mathcal{A}_{n,q}$ and all $\mathbf{s} = (s_0, \ldots, s_{n+1}) \in \mathbf{S}_q^{n+2}$ by

$$\mathbf{s}(\mathbf{a}) := (s_0 a_0 s_1^{-1}, s_1 a_1 s_2^{-1}, \ldots, s_n a_n s_{n+1}^{-1}).$$

Then, for any $F \in \mathcal{B}_b^{\mathrm{sym}}(E_n^q)$ we have

$$(1.12) \qquad\qquad \Delta_{n,q}^{\mathbf{b}}(F) = \Delta_{n,q}^{\mathbf{s}(\mathbf{b})}(F).$$

The identity would not hold if $F$ was not a symmetric function. However, as already mentioned, this is not a serious restriction as far as the determination of the $q$-tensor occupation measures is concerned, since the latter are symmetry-invariant.

DEFINITION 1.4. We write $\mathcal{F}_{n,q}$ for the set of orbits for the action of $\mathbf{S}_q^{n+2}$ on $\mathcal{A}_{n,q}$. We write $\mathbf{a} \approx \mathbf{b}$ if $\mathbf{a}$ and $\mathbf{b}$ belong to the same orbit, that is, if there exists $\mathbf{s} \in \mathbf{S}_q^{n+2}$ such that $\mathbf{a} = \mathbf{s}(\mathbf{b})$.

If $\mathbf{a} \in \mathcal{A}_{n,q}$, the orbit of $\mathbf{a}$ under the action of $\mathbf{s} \in \mathbf{S}_q^{n+2}$ is written $\overline{\mathbf{a}}$. We also write $Stab(\mathbf{a})$ for the stabilizer of $\mathbf{a}$ in $\mathbf{S}_q^{n+2}$. According to the class formula, the number of elements in $\overline{\mathbf{a}}$, written $\#(\mathbf{a})$ or $\#(\overline{\mathbf{a}})$, is given by

$$\#(\overline{\mathbf{a}}) = \frac{(q!)^{n+2}}{|Stab(\mathbf{a})|}.$$

Explicit formulas for $\#(\overline{\mathbf{a}})$ and for the various quantities associated to the action of $\mathbf{S}_q^{n+2}$ on $\mathcal{A}_{n,q}$ such as the number of orbits (and much more) will be given later, and, from the combinatorial point of view, form one of the cores of the article.

Notice that, if $\mathbf{a} \approx \mathbf{b}$, $|\mathbf{a}| = |\mathbf{b}|$, so that the two sequences of maps have the same coalescence sequence and coalescence degree. It follows that the notions of coalescence degree and coalescence sequence go over to the set $\mathcal{F}_{n,q}$. In particular, notation such as $|\overline{\mathbf{a}}|$ or $coal(\overline{\mathbf{a}})$ is well defined. In view of (1.12), for any $\mathbf{f} = \overline{\mathbf{a}}$ in $\mathcal{F}_{n,q}$, we set

$$\Delta_{n,q}^{\mathbf{f}} = \Delta_{n,q}^{\mathbf{a}} \quad \text{and} \quad |\mathbf{f}| = |\mathbf{a}|.$$

The subset of $\mathcal{F}_{n,q}$ associated to $\mathcal{A}_{n,q}(\mathbf{p})$ is written $\mathcal{F}_{n,q}(\mathbf{p})$. We are now in position to state one of the main results of this article.

THEOREM 1.5. *For any $1 \leq q \leq N$ we have the polynomial expansion*

$$\mathbb{Q}_{n,q}^N = \gamma_n^{\otimes q} + \sum_{1 \leq k \leq (q-1)(n+1)} \frac{1}{N^k} \partial^k \mathbb{Q}_{n,q}$$

*with the collection of signed measures $\partial^k \mathbb{Q}_{n,q}$ given by the formula*

$$\partial^k \mathbb{Q}_{n,q}(F) = \sum_{\mathbf{r} \leq \mathbf{q} - \mathbf{1} : \|\mathbf{r}\| = k} \sum_{\mathbf{f} \in \mathcal{F}_{n,q}(\mathbf{r})} \frac{s(|\mathbf{f}|, \mathbf{q} - \mathbf{r}) \#(\mathbf{f})}{(\mathbf{q})_{|\mathbf{f}|}} \Delta_{n,q}^{\mathbf{f}}(F)$$

*for any $F \in \mathcal{B}_b^{\mathrm{sym}}(E_n^q)$.*

In Section 2, we will prove the exact polynomial expansion (in the parameter $N^{-1}$) of Theorem 1.5. In Section 3, we will develop combinatorial concepts in order to compute $\#(\mathbf{f})$ for any $\mathbf{f} \in \mathcal{F}_{n,q}$ (Theorem 3.8). In Section 3, we shall also present several important consequences of Theorem 1.5, including explicit descriptions of the first two order terms in the polynomial expansion, and a new extension of the Wick product formula to forests (Theorem 3.14).

Although more difficult to study since their analysis involves tensor product measures on path-spaces, the distributions $\mathbb{P}_{n,q}^N$ have similar expansions (Theorem 4.12). The description of these expansions, and of the corresponding signed measures, is postponed to Section 4.

## 2. Particle measures expansions.

2.1. *A preliminary stochastic tensor product formula.* The link between the two, usual and restricted, tensor products measures, $(\eta_n^N)^{\otimes q}$ and $(\eta_n^N)^{\odot q}$, relies in the end on a simple observation, that will appear to be fundamental for all our forthcoming computations. Throughout this section, integers $N \geq q \geq 1$ and a measurable state space $E$ are fixed once for all.

Consider first the surjection

$$\pi : [q]^{[q]} \times \langle q, N \rangle \longrightarrow [N]^{[q]},$$

$$(s, a) \longmapsto as := a \circ s.$$

LEMMA 2.1. *Let $b \in [N]^{[q]}$; then the cardinal $|\pi^{-1}(\{b\})|$ of $\pi^{-1}(\{b\})$ only depends on the cardinal $|b|$ of the image of $b$. It is given by*

$$|\pi^{-1}(\{b\})| = (N - |b|)_{q-|b|} (q)_{|b|}.$$

SKETCH OF THE PROOF. Let $b \in [N]^{[q]}$. How many possibilities do we have to write $b$ in the form $b = a \circ s$ with $s \in [q]^{[q]}$, $a \in \langle q, N \rangle$? Set $k_1, \ldots, k_{|b|} \in [q]$ such that $b(\{k_1, \ldots, k_{|b|}\}) = Im(b)$. Any sequence (without repetitions) $s(k_1), \ldots, s(k_{|b|})$ determines uniquely $s$ (since $a$ is injective and $b = a \circ s$) so that there are $(q)_{|b|}$ possibilities for $s$. Then $a(s(k_i)) = b(k_i)$ is fixed for any $i \in \{1, \ldots, |b|\}$ so there remain $(N - |b|)_{q-|b|}$ possibilities for $a$.  □

We denote by $m(\mathbf{x})$ the empirical measure associated with an $N$-uple $\mathbf{x} = (x^i)_{1 \leq i \leq N} \in E^N$:

$$m(\mathbf{x}) := \frac{1}{N} \sum_{i=1}^{N} \delta_{x^i}.$$

For any integer $q \leq N$, we also consider the empirical measures on $E^q$ defined by

$$m(\mathbf{x})^{\otimes q} := \frac{1}{N^q} \sum_{a \in [N]^{[q]}} \delta_{(x^{a(1)}, \ldots, x^{a(q)})},$$

$$m(\mathbf{x})^{\odot q} := \frac{1}{(N)_q} \sum_{a \in \langle q, N \rangle} \delta_{(x^{a(1)}, \ldots, x^{a(q)})}.$$

Although the notation is self-explanatory, notice that, in the sequel, we will write simply $\mathbf{x}^a$ for $(x^{a(1)}, \ldots, x^{a(q)})$ and $\delta_{\mathbf{x}^a}$ for $\delta_{(x^{a(1)}, \ldots, x^{a(q)})}$.

COROLLARY 2.2. *We have*

$$m(\mathbf{x})^{\otimes q} = m(\mathbf{x})^{\odot q} D_{L_q^N},$$

*where*

$$L_q^N = \frac{1}{N^q} \sum_{a \in [q]^{[q]}} \frac{(N)_{|a|}}{(q)_{|a|}} a.$$

PROOF. In view of the following identity, that holds for any function $F \in \mathcal{B}_b^{\text{sym}}(E^q)$, with $q \leq N$:

$$\delta_{\mathbf{x}^a} D_b(F) = D_b(F)(x^{a(1)}, \ldots, x^{a(q)}) = F(x^{ab(1)}, \ldots, x^{ab(q)}) = \delta_{\mathbf{x}^{ab}}(F),$$

the corollary follows from Lemma 2.1 and the identity

$$\frac{1}{(N)_q}(N)_p(N-p)_{q-p} = 1 \qquad \text{where } N \geq q \geq p. \qquad \square$$

COROLLARY 2.3. *The tensor product measure $m(\mathbf{x})^{\otimes q}$ has a Laurent expansion:*

$$m(\mathbf{x})^{\otimes q} = m(\mathbf{x})^{\odot q} \cdot \sum_{0 \leq k < q} \frac{1}{N^k} D_{\partial^k L_q},$$

*where*

$$\partial^k L_q = \sum_{q-k \leq p \leq q} s(p, q-k) \frac{1}{(q)_p} \sum_{a \in [q]_p^{[q]}} a.$$

The corollary follows from Corollary 2.2 and from the Stirling formula

$$(N)_p = \sum_{1 \leq k \leq p} s(p, k) N^k.$$

Assuming from now on that $\mathbf{x}$ is generic (i.e., $x_i \neq x_j$ if $i \neq j$), we have the following corollary. Notice that the assumption that $\mathbf{x}$ is generic is here only to ensure equalities in Corollary 2.4, if it does not hold, we could replace the "=" by "≤."

COROLLARY 2.4. *The following formulas hold for the expansion of $m(\mathbf{x})^{\otimes q}$:*

$$N\|m(\mathbf{x})^{\otimes q} - m(\mathbf{x})^{\odot q}\|_{\text{TV}} \xrightarrow[N \to +\infty]{} q(q-1),$$

$$\|m(\mathbf{x})^{\odot q} D_{\partial^k L_q}\|_{\text{TV}} = \sum_{p=q-k}^{q} |s(p, q-k)| S(q, p).$$

The corollary follows from our previous computations, from the definition of the Stirling numbers of the second kind and from the following lemma.

LEMMA 2.5. *Let $Q = \sum_{0 \leq p \leq q} u_p \sum_{a \in [N]^{[q]}} \delta_{\mathbf{x}^a}$, where the $u_p$ are arbitrary real coefficients; then we have*

$$\|Q\|_{\text{TV}} = \sum_{0 \leq p \leq q} |u_p|(N)_p S(q, p).$$

PROOF. Indeed, we have $\|Q\|_{TV} \leq \sum_{0 \leq p \leq q} |u_p|(N)_p S(q,p)$ by direct inspection, since $\|\delta_{\mathbf{x}^a}\|_{TV} = 1$, for all $a \in [q]^{[q]}$. Let us write $sgn(u_k)$ for the sign of $u_k$, and let us introduce the function $\phi \in \mathcal{B}_b(E^q)$ defined by

$$\phi(y_1, \ldots, y_q) = \sum_{0 \leq k \leq q} sgn(u_k) \delta^k_{\lambda(y_1, \ldots, y_q)},$$

where $\lambda(y_1, \ldots, y_q) := |\{y_1, \ldots, y_q\}|$ and $\delta^k_{\lambda(y_1, \ldots, y_q)}$ stands for the Dirac function. Then, $Q(\phi) = \sum_{0 \leq p \leq q} |u_p|(N)_p S(q,p)$ and the lemma follows. $\square$

2.2. *Laurent expansions.* The goal of this section is to prove Theorem 1.5, that is, to derive a Laurent expansion of the measures $\mathbb{Q}^N_{n,q} \in \mathcal{M}(E^q_n)$, $\mathbb{Q}^N_{n,q}(F) := \mathbb{E}((\gamma^N_n)^{\otimes q}(F))$, $F \in \mathcal{B}_b(E^q_n)$, with respect to the population size parameter $N$. The following proposition is fundamental. Recall from Section 1.3 that we write $\mathcal{F}_{n,q}$ for the set of orbits in $\mathcal{A}_{n,q}$ under the action of the permutation group $\mathbf{S}^{n+2}_q$.

PROPOSITION 2.6. *For any integers $q \leq N$, any time parameter $n \in \mathbb{N}$ and any $F \in \mathcal{B}^{sym}_b(E^q_n)$, we have the Laurent expansion*

$$\mathbb{Q}^N_{n,q}(F) = \frac{1}{N^{q(n+1)}} \sum_{\mathbf{f} \in \mathcal{F}_{n,q}} \frac{(\mathbf{N})_{|\mathbf{f}|}}{(\mathbf{q})_{|\mathbf{f}|}} \#(\mathbf{f}) \Delta^{\mathbf{f}}_{n,q}(F)$$

*with the mappings $\Delta_{n,q}$ introduced in (1.11) properly extended to $\mathcal{F}_{n,q}$.*

PROOF. Combining the definition of the particle model with Corollary 2.2, we first find that

$$\mathbb{E}((\gamma^N_n)^{\otimes q}(F)|\xi^{(N)}_{n-1})$$
$$= (\gamma^N_n(1))^q \times \mathbb{E}((\eta^N_n)^{\otimes q}(F)|\xi^{(N)}_{n-1})$$
$$= (\gamma^N_n(1))^q \times \mathbb{E}((\eta^N_n)^{\odot q} D_{L^N_q}(F)|\xi^{(N)}_{n-1})$$
[by (1.6)] $\quad = (\gamma^N_n(1))^q \times \mathbb{E}(\Phi_n(\eta^N_{n-1})^{\otimes q} D_{L^N_q}(F)|\xi^{(N)}_{n-1}).$

Using the fact that, conditionally to $\xi^{(N)}_{n-1}$ and for $F \in \mathcal{B}^{sym}_b(E^q_n)$,

$$\Phi_n(\eta^N_{n-1})^{\otimes q}(F) = \frac{(\eta^N_{n-1})^{\otimes q}(Q^{\otimes q}_n F)}{(\eta^N_{n-1})^{\otimes q}(Q^{\otimes q}_n(1))} = \frac{(\gamma^N_{n-1})^{\otimes q}(Q^{\otimes q}_n F)}{(\gamma^N_{n-1})^{\otimes q}(Q^{\otimes q}_n(1))}$$

and

$$(\gamma^N_{n-1})^{\otimes q}(Q^{\otimes q}_n(1)) = (\gamma^N_{n-1} Q_n(1))^q = (\gamma^N_{n-1}(G_{n-1}))^q = (\gamma^N_n(1))^q,$$

we arrive at

$$\mathbb{E}((\gamma_n^N)^{\otimes q}(F)|\xi_{n-1}^{(N)}) = (\gamma_n^N(1))^q \times \mathbb{E}\left(\frac{(\gamma_{n-1}^N)^{\otimes q}(Q_n^{\otimes q} D_{L_q^N} F)}{(\gamma_{n-1}^N)^{\otimes q}(Q_n^{\otimes q}(1))}\Big|\xi_{n-1}^{(N)}\right)$$
(2.1)
$$= \mathbb{E}((\gamma_{n-1}^N)^{\otimes q}(Q_n^{\otimes q} D_{L_q^N} F)|\xi_{n-1}^{(N)}).$$

Integrating over the past, this yields that

$$\mathbb{E}((\gamma_n^N)^{\otimes q}(F)) = \mathbb{E}((\gamma_{n-1}^N)^{\otimes q}(Q_n^{\otimes q} D_{L_q^N} F)).$$

Using a simple induction, we readily obtain the formulas

$$\mathbb{E}((\gamma_n^N)^{\otimes q}(F)) = \mathbb{E}(((\gamma_0^N)^{\otimes q} Q_1^{\otimes q} D_{L_q^N} \cdots Q_n^{\otimes q} D_{L_q^N})(F))$$

$$= \mathbb{E}((\eta_0^{\otimes q} D_{L_q^N} Q_1^{\otimes q} D_{L_q^N} \cdots Q_n^{\otimes q} D_{L_q^N})(F))$$

$$= \frac{1}{N^{q(n+1)}} \sum_{\mathbf{a} \in \mathcal{A}_{n,q}} \frac{(\mathbf{N})_{|\mathbf{a}|}}{(\mathbf{q})_{|\mathbf{a}|}} \Delta_{n,q}^{\mathbf{a}}(F),$$

where $\mathcal{A}_{n,q}$ is the set of sequences of maps introduced in Section 1.3. In view of (1.12), we know that

$$\Delta_{n,q}^{\mathbf{b}}(F) = \Delta_{n,q}^{\mathbf{s}(\mathbf{b})}(F)$$

for any $\mathbf{s} \in \mathbf{S}_q^{n+2}$ and the proposition follows. $\square$

We are now in position to derive Theorem 1.5. According to Proposition 2.6

$$\mathbb{Q}_{n,q}^N(F) = \frac{1}{N^{q(n+1)}} \sum_{\mathbf{p} \in [1,q]^{n+1}} \frac{(\mathbf{N})_{\mathbf{p}}}{(\mathbf{q})_{\mathbf{p}}} \sum_{\mathbf{f} \in \mathcal{F}_{n,q}, |\mathbf{f}|=\mathbf{p}} \#(\mathbf{f}) \Delta_{n,q}^{\mathbf{f}}(F).$$

Using the Stirling formula ([3], page 48), we notice that

(2.2) $$(\mathbf{N})_{\mathbf{p}} = \sum_{\mathbf{1} \leq \mathbf{l} \leq \mathbf{p}} s(\mathbf{p}, \mathbf{l}) N^{|\mathbf{l}|}.$$

We also have that

$$\mathbb{Q}_{n,q}^N(F) = \sum_{\mathbf{1} \leq \mathbf{l} \leq \mathbf{p} \leq \mathbf{q}} s(\mathbf{p}, \mathbf{l}) \frac{1}{N^{|\mathbf{q}-\mathbf{l}|}} \frac{1}{(\mathbf{q})_{\mathbf{p}}} \sum_{\mathbf{f} \in \mathcal{F}_{n,q}, |\mathbf{f}|=\mathbf{p}} \#(\mathbf{f}) \Delta_{n,q}^{\mathbf{f}}(F)$$

from which we conclude that

$$\mathbb{Q}_{n,q}^N(F) = \sum_{\mathbf{r} \leq \mathbf{q}-\mathbf{1}} \sum_{\mathbf{q}-\mathbf{r} \leq \mathbf{p} \leq \mathbf{q}} s(\mathbf{p}, \mathbf{q}-\mathbf{r}) \frac{1}{N^{|\mathbf{r}|}} \frac{1}{(\mathbf{q})_{\mathbf{p}}} \sum_{\mathbf{f} \in \mathcal{F}_{n,q} |\mathbf{f}|=\mathbf{p}} \#(\mathbf{f}) \Delta_{n,q}^{\mathbf{f}}(F).$$

The leading term of the development corresponds to the case $\mathbf{r} = \mathbf{0}$. Since $\#(\mathbf{f}) = (q!)^{n+1}$ if $|\mathbf{f}| = \mathbf{q}$, this term is equal to $\gamma_n^{\otimes q}$. Observe then that, for $\mathbf{r} \leq \mathbf{q}-\mathbf{1}, \|\mathbf{r}\| \leq (q-1)(n+1)$. The theorem follows.

COROLLARY 2.7. *The following formula holds for the expansion of $\mathbb{Q}_{n,q}^N$:*

$$\lim_{N \longrightarrow \infty} N \|\mathbb{Q}_{n,q}^N - \gamma_n^{\otimes q}\|_{\mathrm{TV}} \leq (n+1)q(q-1) \times \sup_{\mathbf{f} \in \mathcal{F}_{n,q}} \|\Delta_{n,q}^{\mathbf{f}}\|_{\mathrm{TV}}.$$

REMARK 2.8. If $0 \leq G \leq 1$, then the above supremum is bounded by 1. The above corollary is here to show an estimate of the size of the error term we derive from Theorem 1.5.

PROOF. The corollary follows from the Laurent expansion of Theorem 1.5. The term $1/N$ in the expansion is given by sequences $\mathbf{r}$ s.t. $\|\mathbf{r}\| = 1$. But then, we have two cases to consider, $|\mathbf{f}| > \mathbf{q} - \mathbf{r}$ (which happens if and only if $|\mathbf{f}| = \mathbf{q}$) and $|\mathbf{f}| = \mathbf{q} - \mathbf{r}$.

Then, we have that $\#(\mathbf{f}) = (q!)^{n+1}$ if $|\mathbf{f}| = \mathbf{q}$ and besides, $|\{\mathbf{r}, \|\mathbf{r}\| = 1\}| = n+1$, and for any such $\mathbf{r}$, $s(\mathbf{q}, \mathbf{q} - \mathbf{r}) = -\binom{q}{2}$. In conclusion, the case $|\mathbf{f}| = \mathbf{q}$ contributes to $(n+1)q(q-1)/2$ to the asymptotic evaluation of $N\|\mathbb{Q}_{n,q}^N - \gamma_n^{\otimes q}\|_{\mathrm{TV}}$.

Now, let us consider a sequence $\mathbf{r}$ with $\|\mathbf{r}\| = 1$, for example, the one, written $\mathbf{r}_i$ with $r_i = 1$, for a given $0 \leq i \leq n$. In that case, there is a unique $\mathbf{f}$ with $|\mathbf{f}| = \mathbf{q} - \mathbf{r}$, which is the set of all sequences of maps $\mathbf{a}$ in $\mathcal{A}_{n,q}$ with $a_j \in \mathbf{S}_q, j \neq i$ and $|a_i| = q - 1$. In particular $\#(\mathbf{f}) = (q!)^{n+1}\binom{q}{2}$, so that, on the whole, the case $|\mathbf{f}| = \mathbf{r}_i$ contributes to $q(q-1)/2$ to the evaluation. □

To conclude this section, notice that the above functional expansions also apply to the dot-tensor product measures

$$(\gamma_n^N)^{\odot q}(F) = (\gamma_n^N(1))^q \times (\eta_n^N)^{\odot q}(F).$$

More precisely, by definition of the particle model we have that

$$(2.3) \qquad \mathbb{E}((\gamma_n^N)^{\odot q}(F)) = \mathbb{E}((\gamma_{n-1}^N)^{\otimes q} Q_n^{\otimes q}(F)) = \mathbb{Q}_{n-1,q}^N(Q_n^{\otimes q} F).$$

**3. Combinatorial methods for path enumeration.** In the present section, we face the problem of computing the cardinals $\#(\mathbf{f})$ involved in our Laurent expansions. We also derive various identities relevant for the fine asymptotical analysis of mean-field particle interacting models, such as the number of elements in $\mathcal{F}_{n,q}$ with a given coalescence degree. To express them, we have chosen to use the language of botanical and genealogical trees, both for technical reasons (since forests appear to be the most natural parametrization of elements in $\mathcal{F}_{n,q}$), and also for the clarity of the exposition.

Forests, graphs, colored graphs and their combinatorics have appeared recently or in various situations in probability and analysis, often in connection with problems in mathematical physics. A classical example is given by the forest formula of Zimmermann, which encodes the renormalization process for the divergent integrals showing up in perturbative quantum field theory. The construction is classical, since it goes back to the foundation of modern particle physics (see, e.g., [2]), but has been revisited recently by Connes and Kreimer, who showed that the algebraic study of Feynman graphs (colored by suitable terms appearing in the Lagrangian of the theory) and of associated trees and forests could give rise to a totally new Hopf algebraic approach to the theory [4]. Another example is provided by the Butcher group and Runge–Kutta methods in numerical analysis. Here also, once again, the fine study of the underlying combinatorial and algebraic structures on trees has given rise to a renewal of the subject [1]. Let us mention at last the appearance of the combinatorics of colored graphs in Gaussian matrix integral models. We refer to [7] and subsequent papers by these authors for further references on the subject.

The approach to Feynman–Kac particle models by means of trees, forests and analogous objects which is introduced in the present article is original, to our best knowledge, and the methods developed in this setting seem complementary to the graph and tree-theoretical results that have been obtained in other application fields.

### 3.1. *Jungles, trees and forests.*

#### 3.1.1. *Definitions.*

In this section, we detail the vocabulary of trees, that will be of constant use later in the article. The following definitions are motivated by the parametrization of tensor product measures expansions introduced in the previous sections.

Let $\mathbf{a} = (a_0, a_1, \ldots, a_n)$ be any sequence of maps $a_k : [p_{k+1}] \to [p_k]$, with $k = 0, \ldots, n$. By associating to each element $i$ in $[p_k]$ the point $(i, k)$ in $\mathbb{R}^2$ and for $k > 0$, an oriented arrow from $(i, k)$ to $(a_{k-1}(i), k-1)$, the sequence can be represented graphically by a graph $J(\mathbf{a})$ in the plane, where the edges are allowed to cross.

More generally, let us consider the set $\mathbb{P}$ of graphs in the plane defined by a (possibly empty) finite set $S$ of points in $\mathbb{R}^2$ (the vertices) and a finite subset $A$ of $S^2$ (the arrows or oriented edges) with the following properties: (1) The elements of $S$ are of the form $(x, k), k \in \mathbb{N}$ (notice that, in particular, elements in $S_k := \{s \in S | \exists x, s = (x, k)\}$ are naturally left-to-right ordered). (2) For each $s \in S$, there is at most one arrow starting from $s$. (3) The elements of $A$ are of the form $((x, k+1), (x', k))$, $k \geq 0$. (4) For any $s \in S$, $s = (x, k)$, there is always a (necessarily unique) sequence of arrows joining $s$ to a vertex $(x', 0)$.

Two elements $(S', A')$ and $(S, A)$ of $\mathbb{P}$ are said to be *strongly equivalent* if and only if there exists a bijection $\phi$ between $S$ and $S'$ with the properties: (1) $\phi$ maps bijectively $S_k$ to $S'_k$. (2) $\phi$ induces a bijection between $A$ and $A'$. They are said to be *LR-equivalent* if they are strongly equivalent and moreover: (3) $\phi$ respects the left-to-right ordering.

DEFINITION 3.1. *A jungle is an equivalence class in $\mathbb{P}$ for the LR-equivalence relation. The set of jungles is written $\mathcal{J}$.*

Notice that the process that is described above and maps a sequence of maps **a** to the corresponding jungle [the equivalence class of $J(\mathbf{a})$ in $\mathbb{P}$, still written $J(\mathbf{a})$ by a slight abuse of notation] can be inverted. That is, jungles are actually in bijection with sequences of maps; however, the introduction of jungles is convenient for our purposes, as it will appear below. If **j** is a jungle, we write $Seq(\mathbf{j})$ for the corresponding sequence of maps [so that $\mathbf{j} = J(Seq(\mathbf{j}))$]. In the example in Figure 1, we have $Seq(\mathbf{j}) = (a_0, a_1)$ with $a_0$ (resp. $a_1$) a map from [7] to [3] (resp. [6] to [7]) and, for example, $a_1(5) = 7$ and $a_0(3) = 2$.

Notice that, since the LR-equivalence relation is weaker than the strong equivalence relation, the strong equivalence relation goes over from graphs in $\mathbb{P}$ to jungles.

DEFINITION 3.2. *A forest is an equivalence class of jungles for the strong equivalence relation. The set of all forests is written $\mathcal{F}$. A tree is a connected forest. The set of all trees is written $\mathcal{T}$.*

Here, connected has the usual topological sense (an oriented graph is connected if and only if two arbitrary vertices of the graph can be related by a sequence of adjacent unoriented edges). The forest naturally associated to a sequence of maps **a** and to the corresponding jungle $J(\mathbf{a})$ is written $F(\mathbf{a})$.

A tree can be viewed alternatively as an abstract graph—more precisely as an abstract finite connected oriented graph without loops such that any vertex has at most one outgoing edge. The empty graph $\varnothing$ is viewed as a tree and is called the *empty tree*. The vertices of a tree without incoming edges are called the *leaves*; the vertices with both incoming edges and an outgoing edge are called the *internal vertices*; the (necessarily unique) vertex without outgoing edge is called the *root*. This terminology, as well as other notions introduced below, extend in a straightforward or self-explanatory way from trees to forests, jungles, and so on.

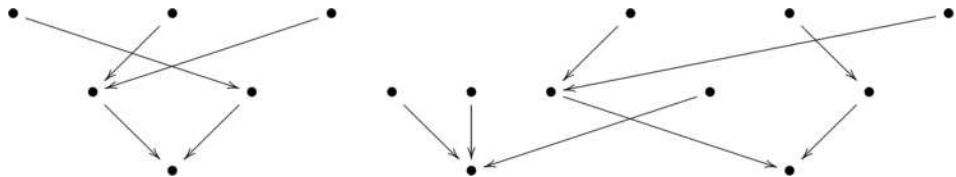

FIG. 1. *The graphical representation of a jungle.*

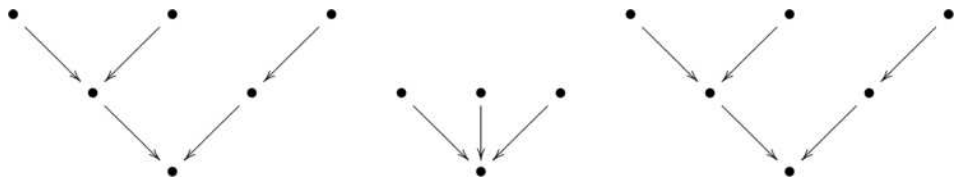

FIG. 2. *The graphical representation of a planar forest* $\mathbf{pf} = PT_1 PT_2 PT_1$.

Forests can be viewed as multisets of trees, that is, sets of trees with repetitions of the same tree allowed. An algebraic notation is convenient for our purposes, and we write

$$\mathbf{f} = T_1^{m_1} \cdots T_k^{m_k}$$

for the forest with the tree $T_i$ appearing with multiplicity $m_i$, $i \leq k$. This algebraic notation will prove useful, among others, when computing Hilbert series for forest enumeration; see Section 3.4. When $T_i \neq T_j$ for $i \neq j$, we say that $\mathbf{f}$ is written in normal form.

We say that a tree $T$ is a subtree of $T'$ and we write $T \subset T'$, if the graph of $T$ is a subgraph of the graph of $T'$, and if the root of $T$ is also the root of $T'$. A subforest $\mathbf{f} \subset \mathbf{f}'$ of a forest $\mathbf{f}'$ is defined accordingly, as a collection of pairwise disjoint subtrees of the trees in $\mathbf{f}'$.

Forests are equivalence classes of jungles. It will be convenient in our forthcoming arguments to choose, for a given forest $\mathbf{f}$, a representative in $\mathcal{J}$ with nice graphical properties. For that purpose, we introduce the notion of *planar forest* and *planar tree*. Namely, a jungle $\mathbf{j}$ such that $Seq(\mathbf{j})$ is a sequence of weakly increasing maps is called a planar forest. A connected planar forest is called a planar tree. Graphically, planar forests are jungles in which the arrows do not cross.

It is important to notice that (as illustrated in Figure 2) a planar forest $\mathbf{pf}$ can be viewed as an ordered sequence of planar trees. Planar forests can therefore be represented by noncommutative monomials (or words) on the set of planar trees. For example, if $PT_1$ and $PT_2$ are two planar trees, $\mathbf{pf} = PT_1 PT_2 PT_1$ is the planar forest obtained by left-to-right concatenation of $PT_1$, $PT_2$ and another copy of $PT_1$.

3.1.2. *Functions and operators on jungles, trees and forests.* The distance between two vertices in a tree is the minimal number of nonoriented edges of a path joining them. The height of a vertex is its distance to the root. We also say that a vertex with height $k$ is a vertex at level $k$ in the tree. The height $ht(T)$ of a tree $T$ is the maximal distance between a leaf and the root. The height $ht(\mathbf{f})$ of a forest $\mathbf{f}$ is the maximal height $ht(T)$ of the trees $T \subset \mathbf{f}$ in the forest. We write $v_k(\mathbf{f})$ for the number of vertices in a forest $\mathbf{f}$ at level $k \geq 0$. The number of vertices at each level in a forest is encoded in the mapping

$$v : \mathbf{f} \in \mathcal{F} \mapsto v(\mathbf{f}) = (v_k(\mathbf{f}))_{k \geq 0}.$$

Notice that for any pair of forests $(\mathbf{f}, \mathbf{g}) \in \mathcal{F}$, we have, using the algebraic notation (the product of forests corresponding to the union of the multisets of trees):

$$v(\mathbf{fg}) = v(\mathbf{f}) + v(\mathbf{g}).$$

We write $\mathbf{V}$ for $v(\mathcal{F})$, which coincides with the set of integer sequences $\mathbf{p} = (p_k)_{k \geq 0} \in \mathbb{N}^{\mathbb{N}}$ satisfying the following property:

$$\exists\, ht(\mathbf{p}) \in \mathbb{N} \text{ s.t.} \quad \inf_{k \leq ht(\mathbf{p})} p_k > 0 \quad \text{and} \quad \sup_{k > ht(\mathbf{p})} p_k = 0.$$

For any $\mathbf{p} \in \mathbf{V}$, we write $\mathcal{F}_{\mathbf{p}}$ for the set of forests $\mathbf{f}$ such that $v(\mathbf{f}) = \mathbf{p}$. Our notation is consistent, since the height $ht(\mathbf{f})$ of a forest $\mathbf{f} \in \mathcal{F}_{\mathbf{p}}$ clearly coincides with the height $ht(\mathbf{p})$ of the integer sequence $\mathbf{p}$. When $\mathbf{p} \in \mathbf{V}$ is chosen so that $p_0 = 1$, the set $\mathcal{F}_{\mathbf{p}}$ reduces to the set $\mathcal{T}_{\mathbf{p}}$ of all trees $T$ such that $v(T) = \mathbf{p}$. These notations on trees and forests go over to planar forests, planar trees and jungles in a self-explanatory way. The notation is also extended to sequences of maps, so that, if $\mathbf{p} = (p_0, \ldots, p_{n+1})$, $\mathcal{A}_{\mathbf{p}}$ stands for the set of sequences $(a_0, \ldots, a_n)$ with $a_i$ a map from $p_{i+1}$ to $p_i$.

For the unit empty tree $T_0 = \varnothing$, we use the conventions

$$v(\varnothing) = \mathbf{0} = (0, 0, \ldots), \qquad ht(\varnothing) = ht(\mathbf{0}) = -1, \qquad \mathcal{F}_{\mathbf{0}} = \{\varnothing\}.$$

For any $n \in (\mathbb{N} \cup \{-1\})$, we denote by $\mathbf{V}_n \subset \mathbf{V}$ the subset of sequences $\mathbf{p}$ such that $ht(\mathbf{p}) = n$. Notice that $\mathbf{V}_{-1} = \{\mathbf{0}\}$. In this notation, the sets of all forests and trees with height $n$ are given by the sets

$$\mathcal{F}_n := \coprod_{\mathbf{p} \in V_n} \mathcal{F}_{\mathbf{p}} \quad \text{and} \quad \mathcal{T}_n := \coprod_{\mathbf{p} \in V_n, p_0 = 1} \mathcal{T}_{\mathbf{p}}.$$

The shift operator

$$B : \mathbf{p} = (p_k)_{k \geq 0} \in \mathbf{V} \mapsto B(\mathbf{p}) = (q_l)_{l \geq 0} \in \mathbf{V} \qquad \text{with } q_l := p_{l+1}$$

induces a canonical bijection, still denoted by $B$, between the set of trees $\mathcal{T}_{\mathbf{p}}$ with $\mathbf{p} \in \mathbf{V}_{n+1}$ and the set of forests $\mathcal{F}_{B(\mathbf{p})}$ obtained by removing the root of the tree:

$$B : T \in \mathcal{T}_{\mathbf{p}} \mapsto B(T) \in \mathcal{F}_{B(\mathbf{p})}.$$

We extend the map $B$ from trees to forests to a map between forests and, if $\mathbf{f} = T_1^{m_1} \cdots T_k^{m_k}$ is a forest, set

$$B(\mathbf{f}) = B(T_1)^{m_1} \cdots B(T_k)^{m_k}.$$

From the point of view of graphs, the operation amounts to removing all the roots and all the edges that have the root as terminal vertex from the graph defining $\mathbf{f}$.

From the graph-theoretic point of view, the number of coalescences $c(\mathbf{f})_k$ at each level $k \geq 0$ in a forest $\mathbf{f} \in \mathcal{F}$ is defined by the mapping

$$c : \mathbf{f} \in \mathcal{F} \mapsto c(\mathbf{f}) := B(v(\mathbf{f})) - |\mathbf{f}|.$$

The sequence $c(\mathbf{f})$ is called the coalescence sequence of $\mathbf{f}$, following the terminology introduced for sequences of maps. In the above displayed formula, the sequence $|\mathbf{f}| = (|\mathbf{f}|_k)_{k \geq 0}$ represents the number of vertices at level $k \geq 0$ minus the number of leaves at the same level, that is, the number of vertices with an ingoing edge, whereas $B(v(\mathbf{f}))_k$ stands, by definition, for the number of edges between the levels $k+1$ and $k$. Finally observe that for any pair of forests $(\mathbf{f}, \mathbf{g}) \in \mathcal{F}$, we have

$$(|\mathbf{fg}| = |\mathbf{f}| + |\mathbf{g}| \text{ and } B(v(\mathbf{f}) + v(\mathbf{g})) = B(v(\mathbf{f})) + B(v(\mathbf{g}))) \Rightarrow c(\mathbf{fg}) = c(\mathbf{f}) + c(\mathbf{g}).$$

We also say that the coalescence order of a vertex in a tree or a forest is its number of incoming edges minus 1. The coalescence degree $\|c(\mathbf{f})\|$ of a tree or a forest $\mathbf{f}$ is the sum of the $c_k(\mathbf{f})$ or, equivalently, the sum of the coalescence orders of its vertices. We say that a tree is *trivial* if its coalescence degree is 0.

3.2. *Automorphism groups on jungles.* For any given sequence of integers $\mathbf{p} \in \mathbf{V}_{n+1}$, the product permutation group

$$\mathbf{S}_\mathbf{p} := (\mathbf{S}_{p_0} \times \mathbf{S}_{p_1} \times \cdots \times \mathbf{S}_{p_{n+1}})$$

acts naturally on sequences of maps $\mathbf{a} = (a_0, a_1, \ldots, a_n) \in \mathcal{A}_\mathbf{p}$, and equivalently on jungles $J(\mathbf{a}) \in \mathcal{J}_\mathbf{p}$ by permutation of the vertices at each level. More formally, for any $\mathbf{s} = (s_0, \ldots, s_{n+1}) \in \mathbf{S}_\mathbf{p}$ this pair of actions is given by

$$\mathbf{s}(\mathbf{a}) := (s_0 a_0 s_1^{-1}, \ldots, s_n a_n s_{n+1}^{-1}) \quad \text{and} \quad \mathbf{s}J(\mathbf{a}) := J(\mathbf{s}(\mathbf{a})).$$

An automorphism $\mathbf{s} \in \mathbf{S}_\mathbf{p}$ of a given jungle $J(\mathbf{a}) \in \mathcal{J}_\mathbf{p}$ is a sequence of permutations that preserves the jungle, in the sense that $\mathbf{s}J(\mathbf{a}) = J(\mathbf{a})$. The set of automorphisms of a given jungle $J(\mathbf{a}) \in \mathcal{J}_\mathbf{p}$ coincides with the stabilizer of $\mathbf{a}$ and $J(\mathbf{a})$ with respect to the group action of $\mathbf{S}_\mathbf{p}$.

By definition of forests, two sequences $\mathbf{a}$ and $\mathbf{b}$ in $\mathcal{A}_\mathbf{p}$ satisfy $F(\mathbf{a}) = F(\mathbf{b})$ if and only if $J(\mathbf{a})$ and $J(\mathbf{b})$ differ only by the order of the vertices, that is, if $\mathbf{a}$ and $\mathbf{b}$ belong to the same orbit under the action of $\mathbf{S}_\mathbf{p}$.

REMARK 3.3. The set of equivalence classes of jungles in $\mathcal{J}_\mathbf{p}$ under the action of the permutation groups $\mathbf{S}_\mathbf{p}$ is in bijection with the set of forests $\mathcal{F}_\mathbf{p}$. In particular, due to the equivalence between the two notions of jungles and sequences of maps, it follows that the set $\mathcal{F}_{n,q}$ of equivalence classes in $\mathcal{A}_{n,q}$ under the natural action of $S_q^{n+2}$ is canonically in bijection with the set of forests $\mathcal{F}_\mathbf{q}$ with $q$ vertices at levels $0 \leq i \leq n+1$.

In the next section, we will compute $\#(\mathbf{f})$ (the number of jungles associated to a given forest $\mathbf{f}$) by computing the cardinal of the stabilizer of a planar forest $\mathbf{pf}$ such that $F(\mathbf{pf}) = \mathbf{f}$ (such a $\mathbf{pf}$ always exists).

EXAMPLE 3.4. The jungle $\mathbf{j}$ represented in Figure 1 is such that $Seq(\mathbf{j}) = (a_0, a_1)$ with

$$a_1(1) = 2, \qquad a_1(2) = a_1(3) = 1, \qquad a_1(4) = a_1(6) = 5, \qquad a_1(5) = 7,$$

$$a_0(1) = a_0(2) = 1, \qquad a_0(3) = a_0(4) = a_0(6) = 2, \qquad a_0(5) = a_0(7) = 3.$$

For convenience, let us represent a permutation in $\mathbf{S}_p$ by the sequence of its values on $1, \ldots, p$, that is, by $(i_1, \ldots, i_q)$ for the permutation $s$ of $[q]$ such that $s(j) = i_j$ (take care that this is not the cycle representation of permutations, in spite of the notational analogy). We get that the jungle represented in Figure 2 is obtained by the action of $(123) \times (1,234,657) \times (312,465)$ on $\mathbf{j}$. The planar forest $\mathbf{j}'$ represented in Figure 2 is thus in the same class as $\mathbf{j}$ (under the action of $\mathbf{S}_{(2,7,6)}$).

A more graphical way of seeing the action of $\mathbf{S}_{(2,7,6)}$ is to say that it permutes the dots while the edges remain attached to their dots.

3.3. *An inductive method for counting jungles.* From the Feynman–Kac mean-field approximation point of view, the elements of $\mathcal{A}_{n,q}$ parametrize the trajectories of families of particles. Two trajectories that are equivalent under the action of $\mathbf{S}_\mathbf{q}$ have the same statistical properties; this is the ground for the formulas for tensor product occupation measures that have been obtained in the previous section.

In the present section, we face the general problem of computing the cardinals $\#(\mathbf{f})$, defined as the number of jungles associated to a given forest $\mathbf{f} \in \mathcal{F}_\mathbf{p}$, for some integer sequence $\mathbf{p} \in \mathbf{V}_{n+1}$. Let us write $\mathbf{f}$ in normal form (with the $T_i$'s all distinct) as a commutative monomial of trees:

$$\mathbf{f} = T_1^{m_1} \cdots T_k^{m_k}, \qquad m_1 + \cdots + m_k = p_0.$$

Let us also choose arbitrary planar trees $PT_i$ representing the $T_i$s, so that we can view $\mathbf{f}$ as the forest associated to the planar forest

$$\mathbf{pf} = (PT_1)^{m_1} \cdots (PT_k)^{m_k}$$

obtained by left-to-right concatenation of $m_1$ copies of $PT_1$, $m_2$ copies of $PT_2, \ldots, m_k$ copies of $PT_k$. We write $PT_1' \cdots PT_{p_0}'$ for the expansion of the noncommutative monomial

(3.1) $$\mathbf{pf} = (PT_1)^{m_1} \cdots (PT_k)^{m_k} = PT_1' \cdots PT_{p_0}'$$

as a product of planar trees (without exponents).

From previous considerations, we have that

$$\#(\mathbf{f}) = \# \, Orb(\mathbf{pf}) = \#\{\mathbf{g} \in \mathcal{J}_\mathbf{p} : \exists \mathbf{s} \in \mathbf{S}_\mathbf{p} \text{ s.t. } \mathbf{g} = \mathbf{s} \cdot \mathbf{pf}\}.$$

Due to the class formula, we also know that

$$\#(\mathbf{f}) = \# \, Orb(\mathbf{pf}) = \frac{|\mathbf{S}_\mathbf{p}|}{Stab(\mathbf{pf})} = \frac{\mathbf{p}!}{Stab(\mathbf{pf})},$$

where $Stab(\mathbf{pf})$ stands for the stabilizer of the jungle $\mathbf{pf} \in \mathcal{J}_\mathbf{p}$.

The computation of $\#(\mathbf{f})$ will be done by induction on $(n+1) = ht(\mathbf{f})$. Let us assume that we know $\#(\mathbf{g})$ for any forest of height less than or equal to $n$. Notice that, then, we also know $\#(\mathbf{t})$ for any tree $\mathbf{t}$ of height $(n+1)$, due to the canonical bijection $B$ between trees and forests: $\#(\mathbf{t}) = \#(B\mathbf{t})$. From the previous discussion, the problem amounts to computing the cardinals of the stabilizers $Stab(\mathbf{pf})$ inductively with respect to the height parameter.

LEMMA 3.5. *There is a natural isomorphism*

$$Stab((PT_1)^{m_1} \cdots (PT_k)^{m_k}) \cap (\{1_{p_0}\} \times \mathbf{S}_{(p_1,\ldots,p_{n+1})})$$
$$\sim Stab(PT'_1) \times \cdots \times Stab(PT'_{p_0})$$

*[with the same notation as in (3.1)].*

The lemma encodes the observation that an automorphism of a planar forest which acts as the identity on the roots decomposes as a product of automorphisms of the planar trees in the forest. This follows, for example, from the more general fact that an automorphism of a jungle necessarily preserves connectedness properties.

Besides, the first coordinate mapping

$$\pi : \mathbf{s} = (s_0, \ldots, s_{n+1}) \in Stab((PT_1)^{m_1} \cdots (PT_k)^{m_k})$$
$$\mapsto \pi(\mathbf{s}) = s_0 \in (\mathbf{S}_{m_1} \times \cdots \times \mathbf{S}_{m_k})$$

is a surjective map (it even has a natural section, the construction of which is omitted), and

$$Ker(\pi) = \pi^{-1}(1_{p_0}) \sim Stab(PT'_1) \times \cdots \times Stab(PT'_{p_0})$$
$$\sim Stab(PT_1)^{m_1} \times \cdots \times Stab(PT_k)^{m_k}.$$

This yields the isomorphism formula

(3.2)
$$Stab((PT_1)^{m_1} \cdots (PT_k)^{m_k})/(Stab(PT_1)^{m_1} \times \cdots \times Stab(PT_k)^{m_k})$$
$$\sim (\mathbf{S}_{m_1} \times \cdots \times \mathbf{S}_{m_k}).$$

Let us introduce some further notation.

DEFINITION 3.6. Let $\mathbf{f}$ be a forest written in normal form $\mathbf{f} = T_1^{m_1} \cdots T_k^{m_k}$. In this situation, we say that the unordered $k$-uplet $(m_1, \ldots, m_k)$ is the symmetry multiset of the tree $T = B^{-1}(\mathbf{f})$ and write

$$\mathbf{s}(T) = \mathbf{s}(B^{-1}(T_1^{m_1} \cdots T_k^{m_k})) = (m_1, \ldots, m_k).$$

The symmetry multiset of the forest $\mathbf{f} = T_1^{m_1} \cdots T_k^{m_k}$ is the disjoint union of the symmetry multisets of its trees

$$\mathbf{s}(T_1^{m_1} \ldots T_k^{m_k}) = \left( \underbrace{\mathbf{s}(T_1), \ldots, \mathbf{s}(T_1)}_{m_1 \text{ terms}}, \ldots, \underbrace{\mathbf{s}(T_k), \ldots, \mathbf{s}(T_k)}_{m_k \text{ terms}} \right).$$

For example, the symmetry multiset of the first tree of the forest $\mathbf{f}$ displayed in Figure 2 is $(1,1)$; the symmetry multiset of $B^{-1}(\mathbf{f})$ is $(2,1)$; and the symmetry multiset of $\mathbf{f}$ is $(1,1,3,1,1)$.

The definition of symmetry multisets extends naturally to planar forests and planar trees. In particular, if $\mathbf{pf}$ is a planar forest with $\mathbf{f}$ as its underlying forest, we set $\mathbf{s}(\mathbf{pf}) := \mathbf{s}(\mathbf{f})$.

We deduce from (3.2) the following recursive formula for the computation of $\#(\mathbf{f})$.

PROPOSITION 3.7. *We have*

$$|Stab((PT_1)^{m_1} \cdots (PT_k)^{m_k})| = \prod_{i=1}^{k} (m_i! |Stab(B(PT_i))|^{m_i})$$

$$= s(B^{-1}(\mathbf{pf}))! \prod_{i=1}^{k} (|Stab(B(PT_i))|^{m_i}).$$

THEOREM 3.8. *The number of jungles (or, equivalently, of Feynman–Kac mean-field type trajectories with the same statistics) in $\mathbf{f} \in \mathcal{F}_{\mathbf{p}}$, with $\mathbf{p} \in \mathbf{V}_{n+1}$, is given by*

$$\#(\mathbf{f}) = \frac{\mathbf{p}!}{\prod_{i=-1}^{n} \mathbf{s}(B^i(\mathbf{f}))!},$$

*where we use the usual multi-index notation to define $\mathbf{s}(B^i(\mathbf{f}))!$.*

PROOF. It clearly suffices to prove that

$$|Stab(\mathbf{pf})| = \prod_{i=-1}^{n} \mathbf{s}(B^i(\mathbf{pf}))!$$

for any planar forest **pf** of height $(n+1)$. We check this assertion by induction on the height parameter. First, we observe that a planar forest **pf** of height 1 can be represented as a noncommutative monomial

$$\mathbf{pf} = (PT_1)^{m_1} \cdots (PT_k)^{m_k}$$

with different planar trees $PT_i$ of height 1, and some sequence of integers $m_i$. In that case, the planar forests $B(PT_i)$ reduce to sequences of trees with null height. This yields that

$$|Stab(B(PT_i))| = \mathbf{s}(PT_i)!.$$

By Proposition 3.7, we conclude that

$$|Stab(\mathbf{pf})| = \prod_{j=1}^{k}(m_j!(\mathbf{s}(PT_j)!)^{m_j}) = \mathbf{s}(B^{-1}(\mathbf{pf}))!\mathbf{s}(\mathbf{pf})!.$$

This ends the proof of the formula at rank 1.

Suppose now that the assertion is satisfied for any planar forest with height at most $n$. By Proposition 3.7, for any planar forest $\mathbf{pf} = (PT_1)^{m_1} \cdots (PT_k)^{m_k}$, with height $(n+1)$, and written in terms of distinct planar trees $PT_i$, we have

$$|Stab(\mathbf{pf})| = \mathbf{s}(B^{-1}(\mathbf{pf}))! \prod_{i=1}^{k} |Stab(B(PT_i))|^{m_i}.$$

Since the planar forests $B(PT_i)$ have height at most $n$, the induction hypothesis implies that

$$|Stab(B(PT_i))| = \prod_{j=-1}^{n-1} \mathbf{s}(B^j(B(PT_i)))! = \prod_{j=0}^{n} \mathbf{s}(B^j(PT_i))!.$$

Recalling that $B^j(\mathbf{pf}) = B^j(PT_1)^{m_1} \cdots B^j(PT_k)^{m_k}$, we also find that

$$\mathbf{s}(B^j(\mathbf{pf})) = (\underbrace{\mathbf{s}(B^j(PT_1)), \ldots, \mathbf{s}(B^j(PT_1))}_{m_1 \text{ terms}}, \ldots, \underbrace{\mathbf{s}(B^j(PT_k)), \ldots, \mathbf{s}(B^j(PT_k))}_{m_k \text{ terms}}).$$

We conclude that

$$|Stab(\mathbf{pf})| = \mathbf{s}(B^{-1}(\mathbf{pf}))! \prod_{j=0}^{n} \mathbf{s}(B^j(\mathbf{pf}))!.$$

The theorem follows. □

3.4. *Hilbert series method for forest enumeration.* In the present subsection, we face the problem of computing other cardinals relevant to the analysis of Feynman–Kac mean-field particle models, according to our Theorem 1.5. For example, we want to be able to compute the number of forests in any $\mathcal{F}_{\mathbf{p}}$ [the set of forests $\mathbf{f}$ such that $v(\mathbf{f}) = \mathbf{p}$], or the number of forests in $\mathcal{F}_{n,q}$ (as defined in Definition 1.4) with a given coalescence degree or a prescribed coalescence sequence. Recall that the notions of coalescence sequence and degree, as defined in the Introduction for sequences of maps in $\mathcal{A}_{n,q}$, go over to arbitrary forests and trees.

Recall first some general properties of free monoids. Let us write $\langle \mathcal{S} \rangle$ for the free commutative monoid generated by a set $\mathcal{S}$:

$$\langle \mathcal{S} \rangle = \{f_1^{m_1} \cdots f_k^{m_k} : k \geq 0 \text{ and } \forall 1 \leq i \leq k, m_i \geq 0 \text{ and } f_i \in \mathcal{S}\}$$

with the convention $f_1^{m_1} \cdots f_k^{m_k} = 1$, for $k = 0$. We then have the Hilbert series (i.e., formal series) expansions

$$(3.3) \qquad \frac{1}{1-s} = \sum_n s^n \quad \text{and} \quad \prod_{s \in \mathcal{S}} \frac{1}{1-s} = \sum_{x \in \langle \mathcal{S} \rangle} x.$$

Now, if $\xi$ is a multiplicative map from $\langle \mathcal{S} \rangle$ to an algebra $A$, we have, whenever the expression on the right-hand side makes sense (this will be obviously the case in the examples that we will consider):

$$\sum_{x \in \langle \mathcal{S} \rangle} \xi(x) = \prod_{s \in \mathcal{S}} \frac{1}{1 - \xi(s)}.$$

Since a forest is, in the algebraic representation, nothing but a commutative monomial of trees, this result applies to forest enumeration. Moreover, as we shall see, suitable refinements of the previous identity lead naturally to a formal series enumeration approach to the various quantities meaningful in the study of Feynman–Kac particle models.

DEFINITION 3.9. The (multidegree) Hilbert series of forests, $\mathbb{H}^n_{\mathcal{F}}(\mathbf{x})$, is the Hilbert series associated to the partition of the set of forests of height less than or equal to $n$ into subsets according to the number of vertices at each level, that is,

$$\mathbb{H}^n_{\mathcal{F}}(\mathbf{x}) := \sum_{ht(\mathbf{p}) \leq n} |\mathcal{F}_{\mathbf{p}}| \mathbf{x}^{\mathbf{p}} = \mathbb{H}^{n-1}_{\mathcal{F}}(\mathbf{x}) + \sum_{\mathbf{p} \in \mathbf{V}_n} |\mathcal{F}_{\mathbf{p}}| \mathbf{x}^{\mathbf{p}},$$

where we write $\mathbf{x}^{\mathbf{p}}$ as a shorthand for $x_0^{p_0} x_1^{p_1} \cdots$.

We let $\overline{\partial}_{\mathbf{p}} = \frac{\partial_{x_0}^{p_0}}{p_0!} \frac{\partial_{x_1}^{p_1}}{p_1!} \cdots$, and we consider the mapping $B^{-1}$ from $\mathbf{V}$ into itself defined for any $k \geq 0$ by

$$B^{-1} : \mathbf{p} \in \mathbf{V}_k \mapsto B^{-1}(\mathbf{p}) = (1, \mathbf{p}) = (1, p_0, p_1, \ldots) \in \mathbf{V}_{k+1}.$$

Since forests $\mathbf{f}$ of height 0 are characterized by the number of roots in $\mathbf{f}$, we clearly have the formula

$$(\forall \mathbf{p} \in \mathbf{V}_0 |\mathcal{F}_\mathbf{p}| = 1) \Longrightarrow \mathbb{H}_\mathcal{F}^0(x_0) = \frac{1}{1 - x_0} = \sum_{p \geq 0} x_0^p.$$

PROPOSITION 3.10.

$$\forall n \geq 1 \qquad \mathbb{H}_\mathcal{F}^n(\mathbf{x}) = \prod_{ht(\mathbf{p}) \leq n-1} \left(\frac{1}{1 - \mathbf{x}^{B^{-1}(\mathbf{p})}}\right)^{|\mathcal{F}_\mathbf{p}|}$$

$$= \mathbb{H}_\mathcal{F}^{n-1}(\mathbf{x}) \times \prod_{\mathbf{p} \in \mathbf{V}_{n-1}} \left(\frac{1}{1 - \mathbf{x}^{B^{-1}(\mathbf{p})}}\right)^{\overline{\partial}_\mathbf{p} \mathbb{H}_\mathcal{F}^{n-1}(\mathbf{0})}.$$

*In addition, the generating function associated with the number of forests with multidegree $\mathbf{p} \in \mathbf{V}_n$, and built with trees with height $n$ (elements of $\mathcal{T}_n$), is given by*

$$\forall n \geq 1 \qquad \sum_{\mathbf{p} \in \mathbf{V}_n} |\langle \mathcal{T}_n \rangle \cap \mathcal{F}_\mathbf{p}| \mathbf{x}^\mathbf{p} = \prod_{\mathbf{p} \in \mathbf{V}_{n-1}} \left(\frac{1}{1 - \mathbf{x}^{B^{-1}(\mathbf{p})}}\right)^{\overline{\partial}_\mathbf{p} \mathbb{H}_\mathcal{F}^{n-1}(\mathbf{0})}.$$

The proposition follows directly from the formula for Hilbert series of free commutative monoids and from the following three observations. First, the map that sends a forest in $\mathcal{F}_\mathbf{p}$ to $\mathbf{x}^\mathbf{p}$ is a multiplicative map from the set of forests (viewed as a free commutative monoid) to the algebra of formal power series over the variables $x_i$. Second, as we have already observed, the set of forests identifies naturally with the free commutative monoid over the set of trees, another formulation of the fact that a forest is a multiset of trees. Third, the maps $B^{-1}$ and $B$ define bijections between the set of forests of height $n - 1$ and trees of height $n$.

Notice that these formulas make the Hilbert series computable at any finite vertex order, and at any height using any formal computation software. The first two orders can be handily computed. We already know the generating series for $n = 0$. By Definition 3.9 we find

$$\mathbb{H}_\mathcal{F}^1(\mathbf{x}) = \sum_{p_0 \geq 0} x_0^{p_0} + \sum_{\mathbf{p} \in \mathbf{V}_1} |\mathcal{F}_\mathbf{p}| x_0^{p_0} x_1^{p_1}.$$

This readily yields that for any $\mathbf{p} \in \mathbf{V}_1$, we have

$$|\mathcal{F}_\mathbf{p}| = \#\left\{(k_n)_{n \geq 0} \in \mathbb{N}^\mathbb{N} : \sum_{n \geq 0} k_n = p_0 \text{ and } \sum_{n \geq 0} nk_n = p_1\right\}.$$

More generally, recall that for any $m \geq 1$, we have the following formal series expansion

$$\left(\frac{1}{1-u}\right)^m = \sum_{k \geq 0} \frac{(m-1+k)!}{(m-1)!k!} u^k$$

so that

$$\mathbb{H}_{\mathcal{F}}^{n+1}(\mathbf{x}) = \prod_{ht(\mathbf{p}) \leq n} \sum_{k_{\mathbf{p}} \geq 0} \frac{(|\mathcal{F}_{\mathbf{p}}|-1+k_{\mathbf{p}})!}{(|\mathcal{F}_{\mathbf{p}}|-1)!k_{\mathbf{p}}!} \mathbf{x}^{k_{\mathbf{p}} B^{-1}(\mathbf{p})}$$

$$= \sum_{k:\{\mathbf{p}:ht(\mathbf{p}) \leq n\} \to \mathbb{N}} \left(\prod_{ht(\mathbf{p}) \leq n} \frac{(|\mathcal{F}_{\mathbf{p}}|-1+k(\mathbf{p}))!}{(|\mathcal{F}_{\mathbf{p}}|-1)!k(\mathbf{p})!}\right) \mathbf{x}^{\sum_{ht(\mathbf{p}) \leq n} k(\mathbf{p}) B^{-1}(\mathbf{p})}.$$

This implies that for any $\mathbf{p} \in \mathbf{V}_{n+1}$, we have

$$(3.4) \qquad |\mathcal{F}_{\mathbf{p}}| = \sum_{\substack{k:\{\mathbf{q}:ht(\mathbf{q}) \leq n\} \to \mathbb{N} \\ \sum_{ht(\mathbf{q}) \leq n} k(\mathbf{q}) B^{-1}(\mathbf{q}) = \mathbf{p}}} \left(\prod_{ht(\mathbf{q}) \leq n} \frac{(|\mathcal{F}_{\mathbf{q}}|-1+k(\mathbf{q}))!}{(|\mathcal{F}_{\mathbf{q}}|-1)!k(\mathbf{q})!}\right).$$

To compute or estimate the Laurent expansions of the distributions $\mathbb{Q}_{n,q}^N$, we are actually interested in a more precise Hilbert series, namely the one taking into account, besides the multidegrees of forests, their coalescence numbers.

DEFINITION 3.11. We denote by $\mathcal{F}_{\mathbf{p}}[\mathbf{q}] := \mathcal{F}_{\mathbf{p}} \cap c^{-1}(\mathbf{q})$ the set of forests in $\mathcal{F}_{\mathbf{p}}$ with a prescribed coalescence sequence $\mathbf{q}$ (we also use the convention $\mathcal{F}_{\mathbf{0}}[\mathbf{0}] = \{\varnothing\}$).

We let $\mathbb{C}_{\mathcal{F}}^n(\mathbf{x}, \mathbf{y})$ be the (multidegree) Hilbert series of forests associated to the partition of the set of forests of height less than or equal to $n$ with prescribed coalescence sequences, that is,

$$\mathbb{C}_{\mathcal{F}}^n(\mathbf{x}, \mathbf{y}) := \sum_{ht(\mathbf{p}) \leq n} \sum_{\mathbf{q} \in C(\mathbf{p})} |\mathcal{F}_{\mathbf{p}}[\mathbf{q}]| \mathbf{x}^{\mathbf{p}} \mathbf{y}^{\mathbf{q}}.$$

In the above display, $C(\mathbf{p})$ stands for the set of coalescence multidegrees:

$$\forall \mathbf{p} \in \mathbf{V} \qquad C(\mathbf{p}) := c(\mathcal{F}_{\mathbf{p}}) = \{\mathbf{q} \in \mathbb{N}^{\mathbb{N}} : (B(\mathbf{p}) - \mathbf{p})_+ \leq \mathbf{q} \leq (B(\mathbf{p}) - 1)_+\}$$

with $(B(\mathbf{p}) - 1)_+ = ((p_k - 1)_+)_{k \geq 1}$, for any $\mathbf{p} = (p_k)_{k \geq 0} \in \mathbf{V}$.

PROPOSITION 3.12. *The multidegree Hilbert series of coalescent forests* $\mathbb{C}_{\mathcal{F}}^n(\mathbf{x}, \mathbf{y})$ *satisfies the recursive formula*

$$\mathbb{C}_{\mathcal{F}}^n(\mathbf{x}, \mathbf{y}) = \mathbb{C}_{\mathcal{F}}^{n-1}(\mathbf{x}, \mathbf{y}) \prod_{\mathbf{p} \in \mathbf{V}_{n-1}} \prod_{\mathbf{q} \in C(B^{-1}(\mathbf{p}))} \left(\frac{1}{1 - \mathbf{x}^{B^{-1}(\mathbf{p})} \mathbf{y}^{\mathbf{q}}}\right)^{\overline{\partial}_{\mathbf{p}} \overline{\partial}_{B(\mathbf{q})} \mathbb{C}_{\mathcal{F}}^{n-1}(\mathbf{0}, \mathbf{0})}.$$

Notice first that $\mathcal{F}_{\mathbf{p}} = \bigcup_{\mathbf{q} \in C(\mathbf{p})} \mathcal{F}_{\mathbf{p}}[\mathbf{q}]$, and $ht(\mathbf{q}) \leq (ht(\mathbf{p}) - 1)$, for any $\mathbf{q} \in C(\mathbf{p})$. Also observe that

$$\mathcal{F}_{B^{-1}(\mathbf{p})} = \mathcal{T}_{(1,\mathbf{p})} \implies C(B^{-1}(\mathbf{p})) = c(\mathcal{T}_{(1,\mathbf{p})}) = \{p_0 - 1\} \times C(\mathbf{p}).$$

Notice at last that the map from forests to formal power series that maps a forest in $\mathcal{F}_{\mathbf{p}}[\mathbf{q}]$ to $\mathbf{x}^{\mathbf{p}}\mathbf{y}^{\mathbf{q}}$ induces a multiplicative map from forests to formal power series. The proposition follows then once again from the formula for Hilbert series of free commutative monoids and the correspondence between trees and forests.

Using the same lines of arguments as before, this proposition can be used to derive a recursive formula for the explicit combinatorial calculation of the number of forests with prescribed heights and coalescence multi-indices. The first two orders are given by

$$\mathbb{C}^0_{\mathcal{F}}(\mathbf{x},\mathbf{y}) = \frac{1}{1 - x_0} \quad \text{and} \quad \mathbb{C}^1_{\mathcal{F}}(\mathbf{x},\mathbf{y}) = \prod_{p_0 \geq 0} \frac{1}{(1 - x_0 x_1^{p_0} y_0^{(p_0-1)_+})}.$$

Notice that the Hilbert series technique can be developed to any order of refinement. For example, it could be used to take into account, besides the number of vertices or of coalescences at each level, the cardinals $\#(\mathbf{f})$, so that the coefficients of the Laurent expansion of the measures $\mathbb{Q}^N_{n,q}$ could be read, in the end, entirely on the corresponding Hilbert series. We leave the task of expressing the recursive formula to the interested reader, and simply point out that the technique allows an easy, systematic, recursive computation of the coefficients of the expansion of the measures $\mathbb{Q}^N_{n,q}$ at any order, in both $n$ and $q$. The observation can be useful, especially in view of the systematic development of numerical schemes and numerical approximations based on Feynman–Kac particle models.

3.5. *Some forests expansions.* In the present section, and the forthcoming one, we take advantage of the language of trees and of the results obtained on their statistics to compute the first orders of the Laurent functional representation of $\mathbb{Q}^N_{n,q}$ and, respectively, to derive a natural generalization to Feynman–Kac particle models of the classical Wick product formula.

Let us fix $n$ and $q \geq 4$, so that the notation $\mathbf{q}$ denotes, once again, the constant sequence of length $n$ associated to $q$. As we have already pointed out, there is only one forest in $\mathcal{F}_{\mathbf{q}}$ without coalescence, which is the product of $q$ trivial trees of height $n$. There is only one forest in $\mathcal{F}_{\mathbf{q}}$ with only one coalescence, at level $i$, that will be written $\mathbf{f}_{1,i}$. Its only nontrivial tree is the tree with one coalescence at level $i$ and two leaves at level $n$. There are two forests with coalescence degree 2 and two coalescences at level $i$, written $\mathbf{f}^1_{2,i}$ and $\mathbf{f}^2_{2,i}$. The notation $\mathbf{f}^1_{2,i}$ denotes the forest with only one nontrivial tree with coalescence degree 2, a vertex with coalescence order 2 at level $i$ and its three leaves at level $n$. The notation $\mathbf{f}_{2,i}^2$ denotes the forest with two nontrivial trees with coalescence degree 1 and the coalescence at level $i$. There are four forests with coalescence degree 2 and two coalescences at levels $i < j$, written $\mathbf{f}_{2,i,j}^k$, $k = 1, \ldots, 4$. The forest $\mathbf{f}_{2,i,j}^1$ has one nontrivial tree with coalescences at levels $i$ and $j$ and its tree leaves at level $n$. The forest $\mathbf{f}_{2,i,j}^2$ has one nontrivial tree with coalescences at levels $i$ and $j$ and its three leaves at the levels $j, n, n$. The forest $\mathbf{f}_{2,i,j}^3$ has two nontrivial trees with one coalescence at level $i$, resp. $j$ and their two leaves at level $n$. The forest $\mathbf{f}_{2,i,j}^4$ has two nontrivial trees with one coalescence at level $i$, resp. $j$ and their two respective leaves at level $j, n$, resp. $n, n$. See Figure 3.

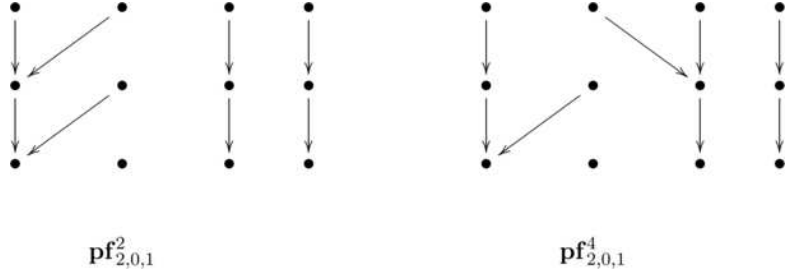

Fig. 3. *Examples of planar tree representatives* $\mathbf{pf}_{r,s}^k$ *of the* $\mathbf{f}_{r,s}^k$'s *for* $q = 4$, $n = 2$.

Expanding the formulas for $\partial^i \mathbb{Q}_{q,n}$, using the formulas obtained in Theorem 3.8 for the cardinals $\#(\mathbf{f})$ and using that $s(q, q-2) = \binom{q}{3}\frac{3q-1}{4}$ (see, e.g., [3], on page 63), we get the following result.

COROLLARY 3.13. *The first three order terms in the polynomial functional representation of* $\mathbb{Q}_{n,q}^N$ *are given by the following formulas:*

$$\partial^0 \mathbb{Q}_{n,q} = \gamma_n^{\otimes q},$$

$$\partial^1 \mathbb{Q}_{n,q} = \frac{q(q-1)}{2} \sum_{0 \leq k \leq n} (\Delta_{n,q}^{\mathbf{f}_{1,k}} - \gamma_n^{\otimes q}),$$

$$\partial^2 \mathbb{Q}_{n,q} = \frac{q!}{(q-3)!3!} \sum_{0 \leq k \leq n} \left( \Delta_{n,q}^{\mathbf{f}_{2,k}^1} + \frac{3}{4}(q-3)\Delta_{n,q}^{\mathbf{f}_{2,k}^2} - \frac{3}{2}(q-1)\Delta_{n,q}^{\mathbf{f}_{1,k}} \right.$$

$$\left. + \frac{(3q-1)}{4}\gamma_n^{\otimes q} \right)$$

$$+ \left( \frac{q(q-1)}{2} \right)^2 \sum_{0 \leq k < l \leq n} (\gamma_n^{\otimes q} - (\Delta_{n,q}^{\mathbf{f}_{1,l}} + \Delta_{n,q}^{\mathbf{f}_{1,k}}))$$

$$+ \frac{q(q-1)}{2} \sum_{0 \leq k < l \leq n} \left( \Delta_{n,q}^{\mathbf{f}_{2,k,l}^2} + (q-2)(\Delta_{n,q}^{\mathbf{f}_{2,k,l}^1} + \Delta_{n,q}^{\mathbf{f}_{2,k,l}^4}) \right.$$

<mark>FEYNMAN–KAC PARTICLE MODELS</mark> 31

$$+ \frac{(q-2)(q-3)}{2} \Delta_{n,q}^{\mathbf{f}_{2,k,l}^3} \Bigg).$$

**3.6. A Wick product formula on forests.** Let now $\mathcal{B}_0^{\mathrm{sym}}(E_n^q) \subset \mathcal{B}_b(E_n^q)$ be the set of symmetric functions $F$ on $E_n^q$ such that

$$(D_{1_{q-1}} \otimes \gamma_n)(F)(x_1,\ldots,x_{q-1}) = \int F(x_1,\ldots,x_{q-1},x_q)\gamma_n(dx_q) = 0.$$

Notice that $\mathcal{B}_0^{\mathrm{sym}}(E_n^q)$ contains the set of functions $F = (P)_{\mathrm{sym}}$, with $P \in \mathrm{Poly}(E_n^q)$, where $\mathrm{Poly}(E_n^q) \subset \mathcal{B}_b(E_n^q)$ stands for the subset of polynomial functions of the form

$$P = \sum_{a \in I} c(a) f^a \qquad \text{with } f^a = (f^{a(1)} \otimes \cdots \otimes f^{a(q)}).$$

In the above display, $I$ is a finite subset of $\mathbb{N}^{[q]}$, $c \in \mathbb{R}^I$, and the elementary functions $f^p \in \mathcal{B}_b(E_n)$ are chosen such that $\gamma_n(f^p) = 0$. For instance, we can take

$$f^p = (g^p - \eta_n(g^p)) \qquad \text{with } g^p \in \mathcal{B}_b(E_n).$$

Notice first that, by Theorem 1.5, for any integer $q$ and any $F \in \mathcal{B}_0^{\mathrm{sym}}(E_0^q)$,

$$\partial^i \mathbb{Q}_{0,q}(F) = 0, \qquad i < \frac{q}{2}$$

and, if $q$ is even,

$$\partial^{q/2} \mathbb{Q}_{0,q}(F) = \frac{q!}{2^{q/2}(q/2)!} \Delta_{0,q}^{\mathbf{f}}(F),$$

where $\mathbf{f}$ is the forest in $\mathcal{F}_{0,q}$ containing $\frac{q}{2}$ copies of the tree of unit height, with two vertices at level 1, and $\frac{q}{2}$ copies of the tree with the root as unique vertex (the reader may check these identities or refer to the more general arguments given below). For symmetric tensor product functions $F = (f^1 \otimes \cdots \otimes f^q)_{\mathrm{sym}}$, associated with a collection of functions $f^i \in \mathcal{B}_b(E_0)$, such that $\eta_0(f^i) = \gamma_0(f^i) = 0$, for any $1 \leq i \leq q$, we readily find that

$$(3.5) \qquad \Delta_{0,q}^{\mathbf{f}}(F) = \frac{2^{q/2}(q/2)!}{q!} \sum_{\mathcal{I}_q} \left( \prod_{\{i,j\} \in \mathcal{I}_q} \eta_0(f^i f^j) \right)$$

as soon as $q$ is even. In the above displayed formula, $\mathcal{I}_q$ ranges over all partitions of $[q]$ into pairs. In a more probabilistic language, the above formula can be interpreted as the $q$th-order central moment

$$(3.6) \qquad \begin{aligned} \partial^{q/2} \mathbb{Q}_{0,q}(F) &= \mathbb{E}(W_0(f^1) \cdots W_0(f^q)) \\ &\stackrel{(\mathrm{Wick})}{=} \sum_{\mathcal{I}_q} \left( \prod_{\{i,j\} \in \mathcal{I}_q} \mathbb{E}(W_0(f^i) W_0(f^j)) \right) \end{aligned}$$

of a Gaussian field $W_0$ on the Banach space of functions $\mathcal{B}_b(E_0)$, such that for any pair of functions $(\varphi, \psi) \in \mathcal{B}_b(E_0)$, $\mathbb{E}(W_0(\varphi)) = 0$ and $\mathbb{E}(W_0(\varphi)W_0(\psi)) = \eta_0(\varphi\psi)$. In the above equation (3.6), the second equality is the classical Wick formula.

We choose to call (3.5) and the more general identities that will appear below "Wick formulas" because of this interpretation. More precisely, in the present section we will show that the Wick formula (3.5) generalizes to forests in $\mathcal{F}_{n,q}$ of arbitrary height. Let us start by listing various straightforward properties of trees and forests. A tree $T$ with coalescence degree $d$ has exactly $(d+1)$ leaves. A forest with coalescence degree $d$ has at most $d$ nontrivial trees, and the equality holds if and only if all its nontrivial trees have coalescence degree 1. In particular, if a forest in $\mathcal{F}_{n,q}$ has coalescence degree $d$, it has at most $(2d)$ leaves belonging to nontrivial trees so that, if $d < \frac{q}{2}$, there is at least one vertex at level $n+1$ belonging to a trivial tree (i.e., a tree with coalescence number 0).

The same reasoning shows that, when $d = \frac{q}{2}$, a forest in $\mathcal{F}_{n,q}$ with coalescence degree $d$, and coalescence sequence $\mathbf{r}$, does not contain a trivial tree of height $(n+1)$ if and only if it is the forest $\mathbf{f_r} := T_0^{r_0} U_0^{r_0} \cdots T_n^{r_n} U_n^{r_n}$, where we write $T_k$ for the unique tree of coalescence degree 1 with a coalescence at level $k$ and its two leaves at level $(n+1)$, and where we write $U_k$ for the trivial tree of height $k$.

We conclude this series of remarks by noting that, if $\mathbf{f} \in \mathcal{F}_{n,q}$ can be written as the product (or disjoint union) of a forest $\mathbf{g}$ in $\mathcal{F}_{n,q-1}$ with $U_{n+1}$, the trivial tree of height $(n+1)$, then, for any $F \in \mathcal{B}_b^{\mathrm{sym}}(E_n^q)$ we have, by definition of the measures $\Delta^{\mathbf{f}}$

$$\Delta_{n,q}^{\mathbf{f}}(F) = \Delta_{n,q}^{\mathbf{g}}(D_{1_{q-1}} \otimes \gamma_n)(F)$$

which is equal to 0 if $F \in \mathcal{B}_0^{\mathrm{sym}}(E_n^q)$.

We are now in position to derive the forest Wick formula using Theorem 1.5.

THEOREM 3.14. *For any even integer $q \leq N$ and any symmetric function $F \in \mathcal{B}_0^{\mathrm{sym}}(E_n^q)$, we have*

$$\forall k < q/2 \quad \partial^k \mathbb{Q}_{n,q}(F) = 0 \quad and$$

(3.7)

$$\partial^{q/2} \mathbb{Q}_{n,q}(F) = \sum_{\mathbf{r} < \mathbf{q}, \|\mathbf{r}\| = \frac{q}{2}} \frac{q!}{2^{q/2} \mathbf{r}!} \Delta_{n,q}^{\mathbf{f_r}} F.$$

*For odd integers $q \leq N$, the partial measure-valued derivatives $\partial^k$ are the null measure on $\mathcal{B}_0^{\mathrm{sym}}(E_n^q)$, up to any order $k \leq \lfloor q/2 \rfloor$.*

We close this section with a Gaussian field interpretation of the Wick formula (3.7). We further assume that $q$ is an even integer. We consider a collection of independent Gaussian fields $(W_k)_{0 \leq k \leq n}$ on the Banach spaces $(\mathcal{B}_b(E_k))_{0 \leq k \leq n}$, with for any $(\varphi_k, \psi_k) \in \mathcal{B}_b(E_k)^2$, and $0 \leq k \leq n$

$$\mathbb{E}(W_k(\varphi_k)) = 0 \quad \text{and} \quad \mathbb{E}(W_k(\varphi_k) W_k(\psi_k)) = \gamma_k(\varphi_k \psi_k).$$

We also introduce the centered Gaussian field $V_n$ on $\mathcal{B}_b(E_n)$ defined for any $\varphi \in \mathcal{B}_b(E_n)$ by the following formula:

$$V_n(\varphi) = \sum_{0 \leq k \leq n} \sqrt{\gamma_k(1)} W_k(Q_{k,n}(\varphi)).$$

Let $(\varphi_i)_{1 \leq i \leq q} \in \mathcal{B}_b(E_n)^q$ be a collection of functions such that $\gamma_n(\varphi_i) = 0$, for any $1 \leq i \leq q$. For the tensor product function $F = \frac{1}{q!} \sum_{\sigma \in \mathbf{S}_q} (\varphi_{\sigma(1)} \otimes \cdots \otimes \varphi_{\sigma(q)})$, one can check that

$$\Delta_{n,q}^{\mathbf{f_r}}(F) = \frac{2^{q/2} \mathbf{r}!}{q!} \sum_{I \in \mathcal{I}} \prod_{0 \leq k \leq n} \sum_{J_k \in \mathcal{I}_k} \left\{ \gamma_k(1)^{r_k} \prod_{(i,j) \in J_k} \gamma_k(Q_{k,n}(\varphi_i) Q_{k,n}(\varphi_j)) \right\}.$$

In the above displayed formula, the first sum is over the set $\mathcal{I}$ of all partitions $I = (I_k)_{0 \leq k \leq n}$ of $[q]$ into $(n+1)$ blocks with cardinality $|I_k| = (2r_k)$, and the second sum ranges over the set $\mathcal{I}_k$ of all partitions $J_k$ of the sets $I_k$ into pairs, with $0 \leq k \leq n$. By definition of the Gaussian fields $(W_k)_{0 \leq k \leq n}$, and due to the classical Wick formula, we find that

$$\frac{q!}{2^{q/2} \mathbf{r}!} \Delta_{n,q}^{\mathbf{f_r}}(F_n) = \sum_{I \in \mathcal{I}} \prod_{0 \leq k \leq n} \mathbb{E} \left( \prod_{i \in I_k} \sqrt{\gamma_k(1)} W_k(Q_{k,n}(\varphi_i)) \right)$$

$$= \mathbb{E} \left( \sum_{I \in \mathcal{I}} \prod_{0 \leq k \leq n} \prod_{i \in I_k} \sqrt{\gamma_k(1)} W_k(Q_{k,n}(\varphi_i)) \right)$$

from which we arrive at

$$\sum_{\mathbf{r} < \mathbf{q}, \|\mathbf{r}\| = q/2} \frac{q!}{2^{q/2} \mathbf{r}!} \Delta_{n,q}^{\mathbf{f_r}}(F_n) = \mathbb{E} \left( \sum_{\mathbf{r}: \|2\mathbf{r}\| = q} \sum_{I} \prod_{0 \leq k \leq n} \prod_{i \in I_k} \sqrt{\gamma_k(1)} W_k(Q_{k,n}(\varphi_i)) \right).$$

Recalling that all Gaussian fields $(W_k)_{0 \leq k \leq n}$ are independent and centered, we get

$$\partial^{q/2} \mathbb{Q}_{n,q}(F_n) = \mathbb{E} \left( \prod_{1 \leq i \leq q} \left( \sum_{0 \leq k \leq n} \sqrt{\gamma_k(1)} W_k(Q_{k,n}(\varphi_i)) \right) \right) = \mathbb{E} \left( \prod_{1 \leq i \leq q} V_n(\varphi_i) \right).$$

Written in a more synthetic way, we have proved the following formula:

$$\partial^{q/2} \mathbb{Q}_{n,q}(F_n) = \mathbb{E}(V_n^{\otimes q}(F_n)).$$

This result can alternatively be derived combining the $\mathbb{L}_q$-mean error estimates presented in [5], Theorem 7.4.2, with the multidimensional central limit theorems presented in [5], Proposition 9.4.1. More precisely, the $q$-dimensional particle random fields

$$(V_n^N(\varphi_n^i))_{1\leq i\leq q} := (\sqrt{N}\gamma_n^N(\varphi_n^i))_{1\leq i\leq q}$$

converge in law to $(V_n(\varphi_n^i))_{1\leq i\leq q}$. By the continuous mapping theorem, combined with simple uniform integrability arguments, one checks that

$$N^{q/2}\mathbb{Q}_{n,q}^N(F_n) = \mathbb{E}((V_n^N)^{\otimes q}(F_n)) \stackrel{N\uparrow\infty}{\longrightarrow} \partial^{q/2}\mathbb{Q}_{n,q}(F_n) = \mathbb{E}(V_n^{\otimes q}(F_n)).$$

**4. Extension to path-space models.** In the present section, we extend our previous analysis to the statistical study of path-spaces. Our goal is now to prove Theorem 4.12 which states an expansion for measures $\mathbb{P}_{n,q}^N$ (analog to Theorem 1.5 for measures $\mathbb{Q}_{n,q}^N$). It will appear in the proof of Theorem 4.12 that we need expansions of tensor products of measures (of the type $\gamma_1^{\otimes k_1} \otimes \cdots \otimes \gamma_n^{\otimes k_n}$). To this purpose, we will introduce and study (in Sections 4.1–4.3) colored trees and forests. We will also define Feynman–Kac semigroup on tensor products of measures (4.4) and finally derive expansions of tensor products of measures (in Section 4.5). Unfortunately, Theorem 4.12 can only be written down after all these technicalities.

Let us make two general remarks about colored forests. First, the (two) colors are used to distinguish between vertices that are authorized or not to have an incoming edge. Second, the statistical meaning of colored forests is similar to the one of forests as far as the approximation measures $\gamma_n^N$ are concerned. Namely, the introduction of colors on vertices is a natural consequence of the replacement of $\gamma_n^N$ by the path-space measures $(\gamma_0^N)^{\otimes q_0} \otimes \cdots \otimes (\gamma_n^N)^{\otimes q_n}$. However, a precise understanding of the meaning of colors can hardly be obtained without entering the details of the constructions, and we refer therefore to Section 4.5 and Theorem 4.2 for deeper insights.

4.1. *Colored trees, forests and jungles.*

DEFINITION 4.1. A colored tree is a tree in $\mathcal{T}$ with colored vertices, with two distinguished colors, say black and white. Only black vertices may have an ingoing edge. That is, equivalently, all white vertices are leaves—the converse being not true in general. A colored forest is a multiset of colored trees. The sets of colored trees and colored forests are denoted, respectively, by $\overline{\mathcal{T}}$ and $\overline{\mathcal{F}}$.

As usual, a colored forest can also be viewed as a commutative monomial over the set of colored trees. In this interpretation, the generating series techniques that we have developed to deal with the enumeration of forests will apply to colored forests. The computation of the corresponding series is left to the interested reader.

Most of the notions associated to trees and forests go over in a straightforward way to colored forests and colored trees. As a general rule, we will write a line over symbols associated to colored trees and colored forests. For instance, we write $\overline{v}_k(\overline{T}) = (v_k(\overline{T}), v'_k(\overline{T}))$ for the number $v_k(\overline{T})$, respectively $v'_k(\overline{T})$, of white, respectively black, vertices in the colored tree $\overline{T}$, at each level $k \geq 0$.

We also let $\overline{\mathbf{V}}$ be the set of all sequences of pairs of integers $\overline{\mathbf{p}} = (\overline{p}_k)_{k \geq 0} \in (\mathbb{N}^2)^{\mathbb{N}}$, with $\overline{p}_k = (p_k, p'_k)$ for every $k \geq 0$, satisfying the following property: $\exists \overline{ht}(\overline{\mathbf{p}}) \in \mathbb{N}$ s.t.

$$\inf_{k \leq \overline{ht}(\overline{\mathbf{p}})} p_k + p'_k > 0, \qquad \inf_{k < \overline{ht}(\overline{\mathbf{p}})} p'_k > 0 \quad \text{and} \quad \sup_{k > \overline{ht}(\overline{\mathbf{p}})} p_k + p'_k = 0.$$

For any $n \in \mathbb{N}$, we denote by $\overline{\mathbf{V}}_n \subset \overline{\mathbf{V}}$ the subset of sequences $\overline{\mathbf{p}}$ such that the height of $\overline{\mathbf{p}}$, $\overline{ht}(\overline{\mathbf{p}})$, is equal to $n$. Finally, we let $\overline{\mathcal{T}}_{\overline{\mathbf{p}}}$ be the set of colored trees $\overline{T}$ with $v_k(\overline{T}) = p_k$ white vertices, and $v'_k(\overline{T}) = p'_k$ black vertices, at each level $k \geq 0$. Since $\overline{T}$ is a colored tree, this implies that $p_0 = 0$ and $p'_0 = 1$, except if $\overline{ht}(\overline{\mathbf{p}}) = 0$. In that case, $p_0$ may also be equal to 1 (and then $p'_0 = 0$).

DEFINITION 4.2. For $\overline{\mathbf{p}} \in \overline{\mathbf{V}}_{n+1}$, let us call, by analogy with the uncolored case, colored jungles (resp. colored planar forests) of type $\overline{\mathbf{p}}$ any sequence of maps (resp. of weakly increasing maps) $(\overline{\alpha}_0, \ldots, \overline{\alpha}_n)$, where

$$\overline{\alpha}_k = (\alpha_k, \alpha'_k) \in [p'_k]^{[p_{k+1}]} \times [p'_k]^{[p'_{k+1}]}.$$

The set of colored jungles is written $\overline{\mathcal{J}}$.

The permutation group analysis of jungles and forests can be extended to the colored case. For that purpose, we write

$$\overline{\mathbf{S}}_{(p,p')} := \mathbf{S}_p \times \mathbf{S}_{p'};$$

the group is acting on $[p]$ by $(\sigma, \beta)(i) := \sigma(i)$ and on $[p']$ by $(\sigma, \beta)(j) := \beta(j)$.

We can define colored planar forests in a graphic way by saying they are colored jungles such that, at any level, edges starting from a fixed color (black or white) do not cross each other. It is the case in Figure 4, although it should be noticed that some vertices are crossing each other (so that the notion of planarity for colored forests is not the intuitive one).

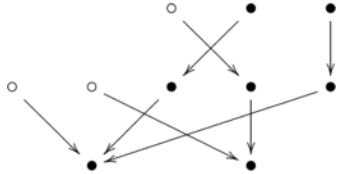

FIG. 4. *A colored jungle which is a colored planar forest.*

Once again, the notation we use on colored forests extends to colored jungles in a self-explanatory way. For instance, $\overline{\mathcal{J}}_{\overline{\mathbf{p}}}$, $\overline{\mathbf{p}} \in \overline{\mathbf{V}}_{n+1}$, stands for the set of colored jungles $\overline{\mathbf{j}}$ with $v_k(\overline{\mathbf{j}}) = p_k$ white vertices, and $v'_k(\overline{\mathbf{j}}) = p'_k$ black vertices, at each level $0 \leq k \leq (n+1)$. For any $\overline{\mathbf{p}} \in \overline{\mathbf{V}}_{n+1}$, we also have that $\overline{\mathcal{J}}_{\overline{\mathbf{p}}} \simeq \overline{\mathcal{A}}_{\overline{\mathbf{p}}}$, where $\overline{\mathcal{A}}_{\overline{\mathbf{p}}} := [p'_0]^{[p_1]} \times [p'_0]^{[p'_1]} \times \cdots \times [p'_n]^{[p_{n+1}]} \times [p'_n]^{[p'_{n+1}]}$.

4.2. *Automorphism groups on colored jungles.* As for the expansion of measures in Section 2, we will need a definition of equivalence classes of colored jungles, a parametrization of the classes and a computation of their cardinals. It appears that most of the constructions in Sections 2 and 3 can be extended without further work to colored jungles. Since all constructions can be mimicked, we only outline the main ideas of the generalization from jungles to colored jungles.

We let $\overline{\mathbf{p}} \in \overline{V}_{n+1}$ be a sequence of $(n+1)$ integer pairs $\overline{p}_k = (p_k, p'_k)$, with $0 \leq k \leq (n+1)$. We associate with $\overline{\mathbf{p}}$ the product permutation group

$$\overline{\mathbf{S}}_{\overline{\mathbf{p}}} = \overline{\mathbf{S}}_{\overline{p}_0} \times \overline{\mathbf{S}}_{\overline{p}_1} \times \cdots \times \overline{\mathbf{S}}_{\overline{p}_{n+1}}.$$

This group acts naturally on sequences of maps $\overline{\mathbf{a}} = (\overline{a}_0, \overline{a}_1, \ldots, \overline{a}_n) \in \overline{\mathcal{A}}_{\overline{\mathbf{p}}}$, where $\overline{a}_i \in [p'_i]^{[p_{i+1}]} \times [p'_i]^{[p'_{i+1}]}$, and on jungles $\overline{J}(\overline{\mathbf{a}})$ in $\overline{\mathcal{J}}_{\overline{\mathbf{p}}}$ by permutation of the colored vertices at each level. More formally, for any $\overline{\mathbf{s}} = (\overline{s}_0, \ldots, \overline{s}_{n+1}) \in \overline{\mathbf{S}}_{\overline{\mathbf{p}}}$ this pair of actions is given by

$$\overline{\mathbf{s}}(\overline{\mathbf{a}}) := (\overline{s}_0 \overline{a}_0 \overline{s}_1^{-1}, \ldots, \overline{s}_n \overline{a}_n \overline{s}_{n+1}^{-1}) \quad \text{and} \quad \overline{\mathbf{s}} \overline{J}(\overline{\mathbf{a}}) := \overline{J}(\overline{\mathbf{s}}(\overline{\mathbf{a}})).$$

Two colored jungles in the same orbit under the action of $\overline{\mathbf{S}}_{\overline{\mathbf{p}}}$ have the same underlying colored forest. Conversely, if two colored jungles in $\overline{\mathcal{J}}_{\overline{\mathbf{p}}}$ have the same underlying colored forests, they differ only by a permutation of the vertices of their colored graphs that preserves the colors of the vertices, and therefore are in the same orbit under the action of $\overline{\mathbf{S}}_{\overline{\mathbf{p}}}$. In other terms, we have the following lemma.

LEMMA 4.3. *Equivalence classes of colored jungles in $\overline{\mathcal{J}}_{\overline{\mathbf{p}}}$ under the action of the permutation groups $\overline{\mathbf{S}}_{\overline{\mathbf{p}}}$ are in bijection with colored forests in $\overline{\mathcal{F}}_{\overline{\mathbf{p}}}$.*

As for usual trees and forests, we write $B$ for the map from colored trees to colored forests defined by removing the root of a colored tree, and we write also, as in Section 3.2, $B$ for the induced map from the set of colored forests into itself.

4.3. *An inductive method for counting colored jungles.* Let us conclude this section by enumerating the number $\#(\overline{\mathbf{f}})$ of colored jungles associated to a given colored forest $\overline{\mathbf{f}} \in \overline{\mathcal{F}_{\overline{\mathbf{p}}}}$, with $\overline{\mathbf{p}} \in \overline{V}_{n+1}$. The process is as in Section 3.2, and the result follows ultimately from the class formula

$$\#(\overline{\mathbf{f}}) = \frac{\overline{\mathbf{p}}!}{|\operatorname{Stab}(\overline{\mathbf{f}})|}.$$

In the above displayed formula, $\overline{\mathbf{p}}!$ stands for the multi-index factorial $\overline{\mathbf{p}}! = \prod_{k=0}^{n+1} p_k! p'_k!$, and we have written abusively $|\operatorname{Stab}(\overline{\mathbf{f}})|$, for the cardinal of the stabilizer in $\overline{\mathbf{S}_{\overline{\mathbf{p}}}}$ of any representative $\overline{\mathbf{pf}}$ of $\overline{\mathbf{f}}$, where $\overline{\mathbf{f}}$ is viewed as an equivalence class of colored jungles.

Let us assume that $\overline{\mathbf{f}}$ can be written, as a monomial over the set of colored trees, as

$$\overline{\mathbf{f}} = \prod_{i=1}^{k} \overline{T_i}^{m_i},$$

where the $\overline{T}_i$'s are pairwise distinct and $m_1 \geq \cdots \geq m_k$. We write $\overline{PT_i}$ for a set of representatives of the $\overline{T_i}$ viewed as equivalence classes of colored jungles (beware that, here, the prime exponent has no color meaning).

As in Section 3.3, we shall write $\mathbf{s}(B^{-1}(\overline{\mathbf{f}}))$ for the unordered $k$-uplet $(m_1, \ldots, m_k)$, called the symmetry multiset of the colored tree $B^{-1}(\overline{\mathbf{f}})$, and extend the notation to colored forests so that $\mathbf{s}(\overline{\mathbf{f}})$ is the disjoint union of the $\mathbf{s}(\overline{T_i})$. The proof in Section 3.3 goes then over without changes, except for the replacement of forests, planar forests and planar trees by their colored analogues. For instance, we have the recursion formula

$$|\operatorname{Stab}(\overline{\mathbf{pf}})| = \prod_{i=1}^{k}(m_i! |\operatorname{Stab} B(\overline{PT_i})|^{m_i}).$$

Specifying the recursion, we deduce the following theorem.

THEOREM 4.4. *The number of colored jungles with underlying colored forest $\overline{\mathbf{f}} \in \overline{\mathcal{F}_{\overline{\mathbf{p}}}}$, with $\overline{\mathbf{p}} \in \overline{V}_{n+1}$, is given by*

$$\#(\overline{\mathbf{f}}) = \frac{\overline{\mathbf{p}}!}{\prod_{i=-1}^{n} \mathbf{s}(B^i(\overline{\mathbf{f}}))!},$$

*where we use the usual multi-index notation to define $\mathbf{s}(B^i(\overline{\mathbf{f}}))!$.*

4.4. *Feynman–Kac semigroups.* In this section, the parameter $n \geq 0$ represents a fixed time horizon. For any sequence of integers $\mathbf{q} = (q_0, \ldots, q_n) \in \mathbb{N}^{n+1}$, and any $-1 \leq m \leq n$, we set

(4.1) $$q'_m := \sum_{m < k \leq n} q_k.$$

Notice that, for $m < n$, $q'_m = q_{m+1} + q'_{m+1}$. We associate with $\mathbf{q}$ the unnormalized Feynman–Kac measures

$$\gamma_n^{\mathbf{q}} = \gamma_0^{\otimes q_0} \otimes \cdots \otimes \gamma_n^{\otimes q_n} \in \mathcal{M}(\mathbf{E}_n^{\mathbf{q}}) \qquad \text{where } \mathbf{E}_n^{\mathbf{q}} = E_0^{q_0} \times \cdots \times E_n^{q_n}$$

is equipped with the tensor product sigma field. Points in $\mathbf{E}_n^{\mathbf{q}}$ are indexed as follows:

$$((x_0^1, \ldots, x_0^{q_0}), \ldots, (x_n^1, \ldots, x_n^{q_n})).$$

DEFINITION 4.5. We let $(\gamma_p^{\mathbf{q}})_{0 \leq p \leq n}$, and $(\mathbf{Q}_p^{\mathbf{q}})_{0 \leq p \leq n}$, be the collection of measures and integral operators defined by

$$\gamma_p^{\mathbf{q}} = (\gamma_0^{\otimes q_0} \otimes \cdots \otimes \gamma_{p-1}^{\otimes q_{p-1}}) \otimes \gamma_p^{\otimes (q_p + q'_p)} \qquad (= \gamma_p^{(q_0, \ldots, q_{p-1}, q_p + q'_p)}),$$

$$\mathbf{Q}_p^{\mathbf{q}} = (D_{1_{q_0}} \otimes \cdots \otimes D_{1_{q_{p-2}}}) \otimes Q_p^{(q_{p-1}, q'_{p-1})},$$

with the operators $Q_p^{(q_{p-1}, q'_{p-1})}$ from $E_{p-1}^{q_{p-1}+q'_{p-1}}$ into $(E_{p-1}^{q_{p-1}} \times E_p^{q'_{p-1}})$ defined by the tensor product formula

$$Q_p^{(q_{p-1}, q'_{p-1})} = D_{1_{q_{p-1}}} \otimes Q_p^{\otimes q'_{p-1}}$$

(recall that the $\gamma_p$'s and $Q_p$'s have been defined in Section 1).

Notice that, for any $0 \leq p \leq n$, $\gamma_p^{\mathbf{q}}$ is a positive measure on the product space

$$\mathbf{E}_p^{\mathbf{q}} := E_0^{q_0} \times \cdots \times E_{p-1}^{q_{p-1}} \times E_p^{q_p + q'_p}$$

and $\mathbf{Q}_p^{\mathbf{q}}$ is a positive integral operator from $\mathbf{E}_{p-1}^{\mathbf{q}}$ into $\mathbf{E}_p^{\mathbf{q}}$.

LEMMA 4.6. For any sequence of integers $\mathbf{q} \in \mathbb{N}^{n+1}$, and any $0 \leq p \leq n$, we have

$$\gamma_n^{\mathbf{q}} = \gamma_p^{\mathbf{q}} \mathbf{Q}_{p,n}^{\mathbf{q}}$$

with the semigroup $(\mathbf{Q}_{p_1, p_2}^{\mathbf{q}})_{1 \leq p_1 \leq p_2 \leq n}$ defined by

$$\mathbf{Q}_{p_1, p_2}^{\mathbf{q}} = \mathbf{Q}_{p_1+1}^{\mathbf{q}} (\mathbf{Q}_{p_1+1, p_2}^{\mathbf{q}}) \qquad (= \mathbf{Q}_{p_1+1}^{\mathbf{q}} \cdots \mathbf{Q}_{p_2-1}^{\mathbf{q}} \mathbf{Q}_{p_2}^{\mathbf{q}})$$

and the convention $\mathbf{Q}_{p_1, p_1}^{\mathbf{q}} = \mathrm{Id}$, the identity operator, for $p_1 = p_2$.

PROOF. We start with observing that

$$\gamma_n^{\mathbf{q}} = \gamma_0^{\otimes q_0} \otimes \gamma_1^{\otimes q_1} \otimes \cdots \otimes \gamma_{n-1}^{\otimes q_{n-1}} \otimes (\gamma_{n-1}^{\otimes q_n} Q_n^{\otimes q_n})$$

$$= \gamma_0^{\otimes q_0} \otimes \gamma_1^{\otimes q_1} \otimes \cdots \otimes \gamma_{n-2}^{\otimes q_{n-2}} \otimes (\gamma_{n-1}^{\otimes (q_{n-1}+q_n)} (D_{1_{q_{n-1}}} \otimes Q_n^{\otimes q_n})).$$

Therefore, we find that $\gamma_n^{\mathbf{q}} = \gamma_{n-1}^{\mathbf{q}} \mathbf{Q}_n^{\mathbf{q}}$. Using a simple induction, the lemma follows. □

4.5. *Unnormalized particle measures.* In this section, we derive a functional representation and a Laurent expansion of particle tensor product measures on the path-space similar to the ones obtained in Theorem 1.5.

In the further development of this section, the time horizon $n \geq 0$ is a fixed parameter, and we let $\mathbf{q} = (q_0, \ldots, q_n) \in \mathbb{N}^{n+1}$ be a given sequence of $(n+1)$ integers. Before we turn to the study of the particle tensor product measures

$$\gamma_n^{\mathbf{q},N} = (\gamma_0^N)^{\otimes q_0} \otimes (\gamma_1^N)^{\otimes q_1} \otimes \cdots \otimes (\gamma_n^N)^{\otimes q_n}$$

and the associated path-space measures $\overline{\mathbb{Q}}_{n,\mathbf{q}}^N \in \mathcal{M}(\mathbf{E}_n^{\mathbf{q}})$ defined by

$$\overline{\mathbb{Q}}_{n,\mathbf{q}}^N : F \in \mathcal{B}_b(\mathbf{E}_n^{\mathbf{q}}) \longmapsto \overline{\mathbb{Q}}_{n,\mathbf{q}}^N(F) = \mathbb{E}(\gamma_n^{\mathbf{q},N}(F)) \in \mathbb{R},$$

let us introduce some useful definitions. We associate with $\mathbf{q}$ the pair of sequences $\overline{\mathbf{q}}$ and $\mathbf{q}'$ defined by

$$\overline{\mathbf{q}} := (q_k, q_k')_{-1 \leq k \leq n} \quad \text{and} \quad \mathbf{q}' := (q_k')_{-1 \leq k \leq n}$$

with the convention $q_{-1} = 0$, and the integer sequence $\mathbf{q}'$ introduced in (4.1). Notice that $q_n' = 0$, $\|\mathbf{q}\| = q_0 + q_0'$, $\overline{ht}(\overline{\mathbf{q}}) = (n+1)$, and

$$|\mathbf{q}'| = \sum_{-1 \leq k \leq n} q_k' = \sum_{0 \leq k \leq n} (q_k + q_k') = \sum_{0 \leq k \leq n} (k+1) q_k.$$

Notice that $[q_{k-1}']^{[q_{k-1}']} \cong [q_{k-1}']^{[q_k]} \times [q_{k-1}']^{[q_k']}$. We will identify the two sets of maps throughout all the sequel.

DEFINITION 4.7. For $\overline{a}_k \in [q_{k-1}']^{[q_k]} \times [q_{k-1}']^{[q_k']}$, we define the Markov transitions

$$\mathbf{D}_{k,\overline{a}_k}^{\mathbf{q}} = D_{1_{q_0}} \otimes \cdots \otimes D_{1_{q_{k-1}}} \otimes D_{\overline{a}_k}.$$

For $u$ a linear combination $\sum_{i \in I} \alpha_i \overline{a}_k^{(i)}$ of $\overline{a}_k^{(i)} \in [q_{k-1}']^{[q_k]} \times [q_{k-1}']^{[q_k']}$, we extend the definition of $\mathbf{D}$ by linearity and write

$$\mathbf{D}_{k,u}^{\mathbf{q}} = \sum_{i \in I} \alpha_i \mathbf{D}_{k,\overline{a}_k^{(i)}}^{\mathbf{q}}.$$

We let $\Delta_{n,\mathbf{q}}$ be the nonnegative measure-valued functional on $\overline{\mathcal{A}}_{\overline{\mathbf{q}}}$ defined by

$$\Delta_{n,\mathbf{q}} : \mathbf{a}' = (\overline{a}_0, \overline{a}_1, \ldots, \overline{a}_n)$$

$$\longmapsto \Delta_{n,\mathbf{q}}^{\mathbf{a}'} := \eta_0^{\otimes |\mathbf{q}|} \mathbf{D}_{0,\overline{a}_0}^{\mathbf{q}} \mathbf{Q}_1^{\mathbf{q}} \mathbf{D}_{1,\overline{a}_1}^{\mathbf{q}} \cdots \mathbf{Q}_n^{\mathbf{q}} \mathbf{D}_{n,\overline{a}_n}^{\mathbf{q}} \in \mathcal{M}(\mathbf{E}_n^{\mathbf{q}}).$$

DEFINITION 4.8. We write $\mathcal{B}_b^{\text{sym}}(\mathbf{E}_n^{\mathbf{q}})$ for the elements $F$ of $\mathcal{B}_b(\mathbf{E}_n^{\mathbf{q}})$ that are symmetric in the first $q_0$ variables, the next $q_1$ variables, ..., and the last $q_n$ variables. That is, we have that

$$\forall (\sigma_0, \ldots, \sigma_n) \in (\mathbf{S}_{q_0} \times \cdots \times \mathbf{S}_{q_n}) \qquad F = (D_{\sigma_0} \otimes \cdots \otimes D_{\sigma_n}) F.$$

We consider the subset $\mathcal{B}_0^{\mathrm{sym}}(\mathbf{E}_n^{\mathbf{q}}) \subset \mathcal{B}_b^{\mathrm{sym}}(\mathbf{E}_n^{\mathbf{q}})$ of all functions $F$ such that

$$(D_{1_{q_0}} \otimes \cdots \otimes D_{1_{q_{p-1}}} \otimes (D_{1_{q_p-1}} \otimes \gamma_p) \otimes D_{1_{q_{p+1}}} \otimes \cdots \otimes D_{1_{q_n}})F = 0$$

for any $0 \leq p \leq n$.

Notice that $\mathcal{B}_0^{\mathrm{sym}}(\mathbf{E}_n^{\mathbf{q}})$ contains the tensor product

$$(\mathcal{B}_0^{\mathrm{sym}}(E_0^{q_0}) \otimes \cdots \otimes \mathcal{B}_0^{\mathrm{sym}}(E_n^{q_n})) = \{F_0 \otimes \cdots \otimes F_n : \forall 0 \leq p \leq n F_p \in \mathcal{B}_0^{\mathrm{sym}}(E_p^{q_p})\}$$

of the sets $\mathcal{B}_0^{\mathrm{sym}}(E_p^{q_p})$ introduced in the beginning of Section 3.6, with $0 \leq p \leq n$. For the same reasons as in (1.12), we also have the symmetry invariance property

$$\forall F \in \mathcal{B}_b^{\mathrm{sym}}(\mathbf{E}_n^{\mathbf{q}}), \forall \overline{\mathbf{s}} \in \overline{\mathbf{S}}_{\overline{\mathbf{q}}}, \forall \overline{\mathbf{j}} \in \overline{\mathcal{J}}_{\overline{\mathbf{q}}} \qquad \Delta_{n,\mathbf{q}}^{\overline{\mathbf{j}}}(F) = \Delta_{n,\mathbf{q}}^{\overline{\mathbf{s}}(\overline{\mathbf{j}})}(F),$$

so that $\Delta_{n,\mathbf{q}}^{\overline{\mathbf{f}}}$ is well defined for $\overline{\mathbf{f}}$ a colored forest. By construction, we also have that

$$\forall F \in \mathcal{B}_0^{\mathrm{sym}}(\mathbf{E}_n^{\mathbf{q}}) \qquad \Delta_{n,\mathbf{q}}^{\overline{\mathbf{f}}}(F) = 0$$

as soon as the colored forest $\overline{\mathbf{f}}$ contains at least one trivial colored tree with a white leaf.

Recall that a tree or a colored tree is said to be trivial if its coalescence sequence is the null sequence of integers. Notice also that for any $\mathbf{l} < \mathbf{q}'$, the set of colored forests with exactly $l_k (\in [0, q'_{k-1}[)$ coalescent edges at level $k$ is given by

$$\overline{\mathcal{F}}_{\overline{\mathbf{q}}}[\mathbf{l}] := \{\overline{\mathbf{f}} \in \overline{\mathcal{F}}_{\overline{\mathbf{q}}} : |\overline{\mathbf{f}}| = \mathbf{q}' - \mathbf{l}\}.$$

Thus, for any $\mathbf{r} < \mathbf{q}'$, the set of colored forests $\overline{\mathcal{F}}_{\overline{\mathbf{q}}}(\mathbf{r})$ with less than $r_k$ coalescent edges at level $k$ is given by

$$\overline{\mathcal{F}}_{\overline{\mathbf{q}}}(\mathbf{r}) := \bigcup_{\mathbf{l} \leq \mathbf{r}} \overline{\mathcal{F}}_{\overline{\mathbf{q}}}[\mathbf{l}].$$

The coalescence degree of a colored forest $\overline{\mathbf{f}} \in \overline{\mathcal{F}}_{\overline{\mathbf{q}}}[\mathbf{l}]$ is the sum $|\mathbf{q}' - \mathbf{l}|$ of the coalescence orders of its vertices.

Notice that a colored forest with coalescence degree $d$ has at most $d$ nontrivial colored trees. In addition, a colored forest with a coalescence degree $d$ has at most $(2d)$ leaves belonging to nontrivial trees. Since a colored forest in $\overline{\mathcal{F}}_{\overline{\mathbf{q}}}$ has exactly $\|\mathbf{q}\|$ white leaves, if $\|\mathbf{q}\| > (2d)$, then it contains at least one trivial colored tree with a white leaf.

Next, we discuss the situation where $\|\mathbf{q}\|$ is an even integer, and we characterize the subset of forests (in $\overline{\mathcal{F}}_{\overline{\mathbf{q}}}$), with a coalescence degree $d = \|\mathbf{q}\|/2$, without any trivial colored tree with a white leaf. This characterization follows the same lines of arguments as the ones presented on page 808. For any integers $0 \leq k \leq l \leq m \leq n$, we let $\overline{T}_{k,l,m}$ be the unique colored tree, with a single coalescence at level $k$, and white leaves at the levels $(l+1)$ and $(m+1)$. Notice that colored forests with $(2r_k)$ white leaves, no black leaves, and coalescent degree $r_k$, with $r_k$ different pairs of coalescent edges at level $k$, are necessarily of the form

$$\overline{\mathbf{f}}_k^{(\mathbf{t}_k)} = \prod_{k \leq l \leq m \leq n} \overline{T}_{k,l,m}^{t_{k,l,m}}$$

for some families of integers $\mathbf{t}_k = (t_{k,l,m})_{k \leq l \leq m \leq n}$ such that $\|\mathbf{t}_k\| = \sum_{k,l,m} t_{k,l,m} = r_k$. Therefore, a colored forest in $\overline{\mathcal{F}}_{\overline{\mathbf{q}}}$, with coalescence degree $\|\mathbf{q}\|/2$, a coalescence sequence $\mathbf{r}$ such that $\|\mathbf{r}\| = \|\mathbf{q}\|/2$, without a trivial tree with a white leaf, has necessarily the following form:

$$\overline{\mathbf{f}}^{(\mathbf{t})} = \overline{\mathbf{f}}_0^{(\mathbf{t}_0)} U_0^{r_0} \cdots \overline{\mathbf{f}}_n^{(\mathbf{t}_n)} U_n^{r_n}$$

for some sequence of families of integers $\mathbf{t} = (\mathbf{t}_k)_{0 \leq k \leq n}$, such that $\|\mathbf{t}_k\| = r_k$, for any $0 \leq k \leq n$. In the above displayed formula, $U_k$ denotes the unique trivial tree with a single black leaf at level $k$, with $0 \leq k \leq n$. We write $\mathbf{t}!$ for $\prod_{k \leq l \leq m \leq n} t_{k,l,m}!$.

We are now ready to extend the Laurent expansions, and the Wick formula presented in Theorem 3.14, to particle models in path-spaces.

THEOREM 4.9. *For any $\mathbf{q} = (q_0, \ldots, q_n)$, with $\|\mathbf{q}\| \leq N$ and any $F \in \mathcal{B}_b^{\mathrm{sym}}(\mathbf{E}_n^{\mathbf{q}})$, we have the Laurent expansion*

$$(4.2) \qquad \overline{\mathbb{Q}}_{n,\mathbf{q}}^N(F) = \left(\gamma_n^{\mathbf{q}} + \sum_{1 \leq k \leq \|(\mathbf{q}'-\mathbf{1})_+\|} \frac{1}{N^k} \partial^k \overline{\mathbb{Q}}_{n,\mathbf{q}}\right)(F)$$

*with the signed measures $\partial^k \overline{\mathbb{Q}}_{n,\mathbf{q}}$ defined by*

$$\partial^k \overline{\mathbb{Q}}_{n,\mathbf{q}}(F) = \sum_{\mathbf{r} < \mathbf{q}', \|\mathbf{r}\|=k} \sum_{\overline{\mathbf{f}} \in \overline{\mathcal{F}}_{\overline{\mathbf{q}}}(\mathbf{r})} \frac{s(|\overline{\mathbf{f}}|, \mathbf{q}' - \mathbf{r}) \#(\overline{\mathbf{f}})}{(\mathbf{q}')_{|\overline{\mathbf{f}}|}} \Delta_{n,\mathbf{q}}^{\overline{\mathbf{f}}}(F).$$

*In addition, for any even integer $\|\mathbf{q}\| \leq N$, and any symmetric function $F \in \mathcal{B}_0^{\mathrm{sym}}(\mathbf{E}_n^{\mathbf{q}})$, we have*

$$\forall k < \|\mathbf{q}\|/2 \qquad \partial^k \overline{\mathbb{Q}}_{n,\mathbf{q}}(F) = 0 \quad and$$

$$\partial^{\|\mathbf{q}\|/2} \overline{\mathbb{Q}}_{n,\mathbf{q}}(F) = \sum_{\mathbf{r} < \mathbf{q}', \|\mathbf{r}\|=\|\mathbf{q}\|/2} \sum_{\langle \mathbf{t} \rangle = \mathbf{r}} \frac{\mathbf{q}!}{2^{\delta(\mathbf{t})} \mathbf{t}!} \Delta_{n,\mathbf{q}}^{\overline{\mathbf{f}}^{(\mathbf{t})}} F$$

*with the integer sequence $\langle \mathbf{t} \rangle := (\|\mathbf{t}_k\|)_{0 \leq k \leq n}$, and the sum of the diagonal terms $\delta(\mathbf{t}) = \sum_{0 \leq k \leq l \leq n} t_{k,l,l}$. For odd integers $\|\mathbf{q}\| \leq N$, the partial measure-valued derivatives $\partial^k$ are the null measure on $\mathcal{B}_0^{\mathrm{sym}}(\mathbf{E}_n^{\mathbf{q}})$, up to any order $k \leq \lfloor \|\mathbf{q}\|/2 \rfloor$.*

Before getting into the details of the proof, we mention that the Wick formula stated above has a natural interpretation in terms of the Gaussian fields $(V_k)_{0 \leq k \leq n}$ introduced in Section 3.6. Following the discussion given in that section, we get that

$$\partial^{\|\mathbf{q}\|/2} \overline{\mathbb{Q}}_{n,\mathbf{q}}(F) = \mathbb{E}((V_0^{\otimes q_0} \otimes \cdots \otimes V_n^{\otimes q_n})(F))$$

for any tensor product function $F$ of the following form: $F = F_{0,q_0} \otimes \cdots \otimes F_{n,q_n}$, with $F_{k,q_k} = \frac{1}{q_k!} \sum_{\sigma_k \in \mathbf{S}_{q_k}} (\varphi_k^{\sigma_k(1)} \otimes \cdots \otimes \varphi_k^{\sigma_k(q_k)})$, $\varphi_k^i \in \mathcal{B}_b(E_k)$ and $\gamma_k(\varphi_k^i) = 0$ for any $1 \leq i \leq q_k$ and any $0 \leq k \leq n$.

PROOF OF THEOREM 4.9. By definition of the particle model, and arguing as in the proof of (2.1), we find that

$$\mathbb{E}([(\gamma_{n-1}^N)^{\otimes q_{n-1}} \otimes (\gamma_n^N)^{\otimes q_n}](\varphi)|\xi_{n-1}^{(N)})$$
$$= [(\gamma_{n-1}^N)^{\otimes q_{n-1}} \otimes \{(\gamma_{n-1}^N)^{\otimes q_n} Q_n^{\otimes q_n} D_{L_{q_n}^N}\}](\varphi)$$
$$= [(\gamma_{n-1}^N)^{\otimes (q_{n-1}+q_n)}] Q_n^{(q_{n-1},q_n)}(D_{1_{q_{n-1}}} \otimes D_{L_{q_n}^N})(\varphi)$$
$$= [(\gamma_{n-1}^N)^{\otimes (q_{n-1}+q'_{n-1})}] Q_n^{(q_{n-1},q'_{n-1})}(D_{1_{q_{n-1}}} \otimes D_{L_{q'_{n-1}}^N})(\varphi)$$

for any $\varphi \in \mathcal{B}_b^{\mathrm{sym}}(E_{n-1}^{q_{n-1}} \times E_n^{q_n})$. This yields that (recalling Definition 4.7)

$$\mathbb{E}(\gamma_n^{\mathbf{q},N}(F)|\xi_{n-1}^{(N)}) = \gamma_{n-1}^{\mathbf{q},N} \mathbf{Q}_n^{\mathbf{q}} \mathbf{D}_{n,L_{q'_{n-1}}^N}^{\mathbf{q}}(F),$$

from which we readily conclude that

$$\mathbb{E}(\gamma_n^{\mathbf{q},N}(F)) = \mathbb{E}(\gamma_{n-1}^{\mathbf{q},N}[\mathbf{Q}_n^{\mathbf{q}} \mathbf{D}_{n,L_{q'_{n-1}}^N}^{\mathbf{q}}(F)]).$$

A simple induction yields that

$$\mathbb{E}(\gamma_n^{\mathbf{q},N}(F)) = \mathbb{E}(\eta_0^{\otimes q_0 + q'_0} \mathbf{D}_{0,L_{q_0+q'_0}^N}^{\mathbf{q}} \mathbf{Q}_1^{\mathbf{q}} \mathbf{D}_{1,L_{q'_0}^N}^{\mathbf{q}} \cdots \mathbf{Q}_n^{\mathbf{q}} \mathbf{D}_{n,L_{q'_{n-1}}^N}^{\mathbf{q}}(F))$$
$$= \frac{1}{N^{|\mathbf{q}'|}} \sum_{\overline{\mathbf{a}} \in \overline{\mathcal{A}}_{\overline{\mathbf{q}}}} \frac{(\mathbf{N})_{|\overline{\mathbf{a}}|}}{(\mathbf{q}')_{|\overline{\mathbf{a}}|}} \Delta_{n,\mathbf{q}}^{\overline{\mathbf{a}}}(F) = \frac{1}{N^{|\mathbf{q}'|}} \sum_{\overline{\mathbf{f}} \in \overline{\mathcal{F}}_{\overline{\mathbf{q}}}} \frac{(\mathbf{N})_{|\overline{\mathbf{f}}|}}{(\mathbf{q}')_{|\overline{\mathbf{f}}|}} \#(\overline{\mathbf{f}}) \Delta_{n,\mathbf{q}}^{\overline{\mathbf{f}}}(F),$$

from which we find the following formula:

$$\overline{\mathbb{Q}}_{n,\mathbf{q}}^N(F) = \frac{1}{N^{|\mathbf{q}'|}} \sum_{\mathbf{1} \leq \mathbf{p} \leq \mathbf{q}'} \frac{(\mathbf{N})_{|\mathbf{p}|}}{(\mathbf{q}')_{|\mathbf{p}|}} \sum_{\overline{\mathbf{f}} \in \overline{\mathcal{F}}_{\overline{\mathbf{q}}}: |\overline{\mathbf{f}}| = \mathbf{p}} \#(\overline{\mathbf{f}}) \Delta_{n,\mathbf{q}}^{\overline{\mathbf{f}}}(F).$$

Using the Stirling formula (2.2), we readily check that

$$\overline{\mathbb{Q}}_{n,\mathbf{q}}^N = \sum_{\mathbf{1} \leq \mathbf{l} \leq \mathbf{q}'} \sum_{\mathbf{1} \leq \mathbf{p} \leq \mathbf{q}'} s(\mathbf{p}, \mathbf{l}) \frac{1}{N^{|\mathbf{q}'-\mathbf{l}|}} \frac{1}{(\mathbf{q}')_{|\mathbf{p}|}} \sum_{\overline{\mathbf{f}} \in \overline{\mathcal{F}}_{\overline{\mathbf{q}}}: |\overline{\mathbf{f}}| = \mathbf{p}} \#(\overline{\mathbf{f}}) \Delta_{n,\mathbf{q}}^{\overline{\mathbf{f}}}.$$

From previous computations, we conclude that

$$\overline{\mathbb{Q}}_{n,\mathbf{q}}^N = \sum_{\mathbf{r}<\mathbf{q}'} \frac{1}{N^{|\mathbf{r}|}} \sum_{\overline{\mathbf{f}}\in\overline{\mathcal{F}}_{\overline{\mathbf{q}}}(\mathbf{r})} s(|\overline{\mathbf{f}}|, \mathbf{q}' - \mathbf{r}) \frac{1}{(\mathbf{q}')_{|\overline{\mathbf{f}}|}} \#(\overline{\mathbf{f}}) \Delta_{n,\mathbf{q}}^{\overline{\mathbf{f}}}.$$

Finally, we notice that $\overline{\mathcal{F}}_{\overline{\mathbf{q}}}(\mathbf{0})$ reduces to the single class of all sequences of bijections in $\overline{\mathcal{A}}_{\overline{\mathbf{q}}}$. The end of the proof of the first assertion is now clear. To end the proof of the theorem, we notice that [with $\mathbf{t}$ as in the expansion of $\partial^{\|\mathbf{q}\|/2}\overline{\mathbb{Q}}_{n,\mathbf{q}}(F)$]

$$\#(\overline{\mathbf{f}}^{(\mathbf{t})}) = \frac{\overline{\mathbf{q}}!}{\prod_{0\leq k\leq n}[r_k!(\prod_{k\leq l\leq m\leq n} t_{k,l,m}!)(\prod_{k\leq l\leq n} 2^{t_{k,l,l}})]}.$$

Since we have $\mathbf{s}(|\overline{\mathbf{f}}^{(\mathbf{t})}|, \mathbf{q}' - \mathbf{r}) = 1$, and $(\mathbf{q}')_{|\overline{\mathbf{f}}^{(\mathbf{t})}|} = \prod_{-1\leq k<n} (q'_k)_{q'_k - r_{k+1}}$, for any $\mathbf{r} = (r_k)_{0\leq k\leq n}$, we conclude that

$$\frac{\mathbf{s}(|\overline{\mathbf{f}}^{(\mathbf{t})}|, \mathbf{q}' - \mathbf{r})}{(\mathbf{q}')_{|\overline{\mathbf{f}}^{(\mathbf{t})}|}} \#(\overline{\mathbf{f}}^{(\mathbf{t})}) = \frac{\mathbf{q}!}{2^{\delta(\mathbf{t})} \mathbf{t}!}$$

with $\mathbf{t}! = \prod_{0\leq k\leq l\leq m\leq n} t_{k,l,m}!$. The end of the proof of the theorem is now straightforward. □

4.6. *Propagations of chaos-type expansions.* This section is essentially concerned with applications of the differential forest expansion machinery developed earlier to propagations of chaos properties of interacting particle models. In order to state and prove the main results of this section, we need to introduce some notation. We shall work throughout with a fixed time horizon $n \geq 0$, and a constant particle block size $q \leq N$. We let $\mathbb{N}_q^{n+1}$ be the set of integer sequences $\mathbf{p} = (p_k)_{0\leq k\leq n} \in \mathbb{N}^{n+1}$ such that $\|\mathbf{p}\| = q$. We associate to a given $\mathbf{p} = (p_k)_{0\leq k\leq n} \in \mathbb{N}^{n+1}$ the pair of integer sequences $\mathbf{p}'$ and $\mathbf{p} + q$ defined by

$$\forall 0 \leq k \leq n \quad p'_k := \sum_{k<l\leq n} p_l \quad \text{and} \quad \mathbf{p} + q := (p_0, \ldots, p_{n-1}, p_n + q).$$

We also denote by $(\overline{G}_k)_{0\leq k\leq n}$ and $\overline{G}_n^{\otimes \mathbf{p}}$ the collection of functions defined by

$$\overline{G}_k := \frac{1}{\gamma_k(G_k)} \times (\eta_k(G_k) - G_k) \quad \text{and} \quad \overline{G}_n^{\otimes \mathbf{p}} := \overline{G}_0^{\otimes p_0} \otimes \cdots \otimes \overline{G}_n^{\otimes p_n}$$

(recall the $G_k$'s have been defined in Section 1). Finally, we consider the integral operators $\overline{\mathbf{Q}}_n^{(\mathbf{p}+q)}$ from $\mathbf{E}_n^{\mathbf{p}+q}$ into $E_{n+1}^q$, defined for any function $F_{n+1} \in \mathcal{B}_b(E_{n+1}^q)$ by the following formula

$$\overline{\mathbf{Q}}_n^{(\mathbf{p}+q)}(F_{n+1}) := \overline{G}_0^{\otimes p_0} \otimes \cdots \otimes \overline{G}_{n-1}^{\otimes p_{n-1}} \otimes (\overline{G}_n^{\otimes p_n} \otimes Q_{n+1}^{\otimes q} \overline{F}_{n+1})_{\text{sym}}$$

with

$$\overline{F}_n := \frac{1}{\gamma_n(1)^q}(F_n - \eta_n^{\otimes q}(F_n)).$$

The following proposition and lemma are technical results that will be needed in the proof of Theorem 4.12.

PROPOSITION 4.10. *For any $q \leq N$ and any $n \geq 0$, we have the polynomial decompositions*

$$E_{q,n}^N := \mathbb{E}((1 - \gamma_n^N(G_n)/\gamma_n(G_n))^q) = \sum_{q/2 \leq k \leq (n+1)(q-1)} \frac{1}{N^k} \partial^k E_{q,n}.$$

*The derivatives of order $q/2 \leq k \leq (n+1)(q-1)$ are given by the following formula:*

$$\partial^k E_{q,n} = \sum_{\mathbf{p} \in \mathbb{N}_q^{n+1}} \sum_{k \leq \|(\mathbf{p}'-1)_+\|} \frac{q!}{\mathbf{p}!} \partial^k \overline{\mathbb{Q}}_{n,\mathbf{p}}(\overline{G}_n^{\otimes \mathbf{p}}).$$

PROOF. We first check the following decomposition:

(4.3) $$1 - \gamma_n^N(G_n)/\gamma_n(G_n) = \sum_{0 \leq p \leq n} \gamma_p^N(\overline{G}_p).$$

To prove this formula, we notice that

$$\begin{aligned}
\gamma_n^N(G_n) - \gamma_n(G_n) &= \prod_{0 \leq p \leq n} \eta_p^N(G_p) - \prod_{0 \leq p \leq n} \eta_p(G_p) \\
&= \gamma_n^N(G_n - \eta_n(G_n)) \\
&\quad + [\gamma_{n-1}^N(G_{n-1}) - \gamma_{n-1}(G_{n-1})] \times \eta_n(G_n) = \cdots \\
&= \sum_{0 \leq p \leq n} \gamma_p^N(G_p - \eta_p(G_p)) \times \left[\prod_{p+1 \leq k \leq n} \eta_k(G_k)\right]
\end{aligned}$$

with the convention $\prod_\varnothing = 1$. Finally, combining the multinomial decomposition

(4.4) $$(1 - \gamma_n^N(G_n)/\gamma_n(G_n))^q = \sum_{\|\mathbf{p}\|=q} \frac{q!}{\mathbf{p}!} \gamma_n^{\mathbf{p},N}(\overline{G}_n^{\mathbf{p}})$$

with the Wick expansion stated in Theorem 4.9 we conclude that

$$E_{q,n}^N = \sum_{\mathbf{p} \in \mathbb{N}_q^{n+1}} \frac{q!}{\mathbf{p}!} \overline{\mathbb{Q}}_{n,\mathbf{p}}^N(\overline{G}_n^{\mathbf{p}}) = \sum_{\mathbf{p} \in \mathbb{N}_q^{n+1}} \sum_{q/2 \leq k \leq \|(\mathbf{p}'-1)_+\|} \frac{1}{N^k} \frac{q!}{\mathbf{p}!} \partial^k \overline{\mathbb{Q}}_{n,\mathbf{p}}(\overline{G}_n^{\mathbf{p}}).$$

This completes the proof of the proposition. □

LEMMA 4.11. *For any $u \in \mathbb{R} - \{1\}$, and $m \geq 1$, we have the decomposition*

$$\frac{1}{(1-u)^{q+1}} = \sum_{0 \leq k \leq m} (q+k)_k \frac{u^k}{k!} + u^m \sum_{1 \leq k \leq q+1} \binom{(q+1)+m}{k+m} \left(\left(\frac{u}{1-u}\right)^k\right)$$

*with $(q+k)_k = (q+k)!/q!$.*

PROOF. This formula seems classical but we did not find it in the literature. The proof is essentially based on the fact that for any $n \geq 0$ and $u \neq 1$, we have

$$f(u) = \frac{1}{1-u} = \sum_{0 \leq k \leq n} u^k + \frac{u^{n+1}}{1-u} \quad \text{and} \quad \frac{\partial^n f}{\partial u^n} = n! f^{n+1}.$$

This implies that for any $m > n$, we have

$$\frac{\partial^n f}{\partial u^n} = \frac{n!}{(1-u)^{n+1}} = \sum_{n \leq k \leq m} (k)_n u^{k-n} + \frac{\partial^n}{\partial u^n}\left(\frac{u^{m+1}}{1-u}\right).$$

Applying the Leibniz binomial derivation formula for $\frac{\partial^n}{\partial u^n}(fg)$ to the product $u^{m+1} \times \frac{1}{1-u}$ we find that

$$\frac{1}{n!} \frac{\partial^n}{\partial u^n}\left(\frac{u^{m+1}}{1-u}\right) = u^{m-n} \sum_{0 \leq k \leq n} \binom{m+1}{n-k} \frac{u^{k+1}}{(1-u)^{k+1}}.$$

The lemma follows. □

THEOREM 4.12. *For any $n \geq 0$, the sequence of probability measures $(\mathbb{P}^N_{n+1,q})_{N \geq q}$ is differentiable up to any order with $\partial^0 \mathbb{P}_{n+1,q} = \eta^{\otimes q}_{n+1}$, and the partial derivatives given by the following formula:*

$$\partial^k \mathbb{P}_{n+1,q} = \sum_{\mathbf{p} \in \bigcup_{0 \leq l < 2k} \mathbb{N}^{n+1}_l} \frac{((q-1)+|\mathbf{p}|)!}{(q-1)!\mathbf{p}!} \partial^k \overline{\mathbb{Q}}_{n,\mathbf{p}+q} \overline{\mathbf{Q}}^{(\mathbf{p}+q)}_n.$$

*At any order $N \geq q$, we have the exact formula*

$$\mathbb{P}^N_{n+1,q} = \eta^{\otimes q}_{n+1} + \sum_{1 \leq k < \lfloor (N-q)/2 \rfloor} \frac{1}{N^k} \partial^k \mathbb{P}_{n+1,q} + \mathbb{R}^N_{n+1,q}$$

*with a remainder measure $\mathbb{R}^N_{n+1,q}$ such that $\sup_{N \geq 1} N^{((N-q)+1)/2} \|\mathbb{R}^N_{n+1,q}\|_{\mathrm{TV}} < \infty$.*

REMARK 4.13. Contrary to the expression in Theorem 1.5, where we have an exact expansion as a sum of polynomial terms, the expansion above,

though exact, is the sum of a polynomial term plus a remainder term which is not polynomial. This particular form comes from (4.5) appearing below in the proof of the theorem, where we have to develop a fraction, whose expansion is an infinite sum.

PROOF OF THEOREM 4.12. We let $F_{n+1}$ be a bounded measurable symmetric function on $\mathbf{E}_{n+1}^q$, such that $\eta_{n+1}^{\otimes q}(F_{n+1}) = 0$. By definition of the particle model, we have that

$$
\begin{aligned}
\mathbb{E}((\eta_{n+1}^N)^{\odot q}(F_{n+1})|\xi_n^{(N)}) &= \gamma_n^N(G_n)^{-q} \times (\gamma_n^N)^{\otimes q} Q_{n+1}^{\otimes q}(F_{n+1}) \\
&= (1 - u_n^N)^{-q} \times (\gamma_n^N)^{\otimes q} Q_{n+1}^{\otimes q}(\overline{F}_{n+1})
\end{aligned}
$$
(4.5)

with the sequence of random variables $u_n^N := (1 - \gamma_n^N(G_n)/\gamma_n(G_n))$. Using Lemma 4.11, we find that

$$
\mathbb{E}((\eta_{n+1}^N)^{\odot q}(F_{n+1})|\xi_n^{(N)}) = \sum_{0 \leq k \leq m} ((q-1)+k)_k \frac{1}{k!}\left(1 - \frac{\gamma_n^N(G_n)}{\gamma_n(G_n)}\right)^k
$$
$$
\times (\gamma_n^N)^{\otimes q}(Q_{n+1}^{\otimes q}\overline{F}_{n+1}) + R_{m,n}^{q,N}(F_{n+1})
$$

with the remainder term

$$
R_{m,n}^{q,N}(F_{n+1}) = \left(1 - \frac{\gamma_n^N(G_n)}{\gamma_n(G_n)}\right)^{m+1}(\gamma_n^N)^{\otimes q}(Q_{n+1}^{\otimes q}\overline{F}_{n+1})
$$
$$
\times \sum_{1 \leq k \leq q} \binom{q+m}{k+m}\left(1 - \frac{\gamma_n^N(G_n)}{\gamma_n(G_n)}\right)^{k-1} \bigg/ \left(\frac{\gamma_n^N(G_n)}{\gamma_n(G_n)}\right)^k.
$$

We also have from Proposition 4.10

$$\sup_{N \geq 1} \sqrt{N} \mathbb{E}(|\gamma_n^N(G_n) - \gamma_n(G_n)|^{m+1})^{1/m+1} < \infty.$$

Using the regularity hypothesis (1.1) on the potential functions we conclude that

$$\sup_{N \geq 1} N^{(m+1)/2} \|R_{m,n}^{q,N}\|_{\mathrm{TV}} < \infty.$$

Finally, using the multinomial decomposition (4.4), we conclude that

$$
\sum_{0 \leq k \leq m} ((q-1)+k)_k \frac{1}{k!}\left(1 - \frac{\gamma_n^N(G_n)}{\gamma_n(G_n)}\right)^k (\gamma_n^N)^{\otimes q}(Q_{n+1}^{\otimes q}\overline{F}_{n+1})
$$
$$
= \sum_{0 \leq k \leq m} ((q-1)+k)_k \sum_{\mathbf{p} \in \mathbb{N}_k^{n+1}} \frac{1}{\mathbf{p}!}\gamma_n^{(\mathbf{p}+q),N}(\overline{\mathbf{Q}}_n^{(\mathbf{p}+q)}(F_{n+1})).
$$

This yields, for any $m + q \leq N$, the functional expansion

$$
\mathbb{P}_{n+1,q}^N = \sum_{\|\mathbf{p}\| \leq m} \frac{((q-1)+|\mathbf{p}|)!}{(q-1)!\mathbf{p}!} \overline{\mathbb{Q}}_{n,\mathbf{p}+q}^N \overline{\mathbf{Q}}_n^{(\mathbf{p}+q)} + R_{m,n}^{q,N}
$$

with a remainder measure $R_{m,n}^{q,N}$, such that $\sup_{N\geq 1} N^{(m+1)/2}\|R_{m,n}^{q,N}\|_{\mathrm{TV}} < \infty$. This implies that for any $k < (m+1)/2$ we have

$$\partial^k \mathbb{P}_{n+1,q} = \sum_{\|\mathbf{p}\|\leq m} \frac{((q-1)+|\mathbf{p}|)!}{(q-1)!\mathbf{p}!} \partial^k \overline{\mathbb{Q}}_{n,\mathbf{p}+q} \overline{\mathbf{Q}}_n^{(\mathbf{p}+q)}.$$

This proves the first assertion of the theorem. Notice that $k$th-order derivative measure $\partial^k \overline{\mathbb{Q}}_{n,\mathbf{p}+q}$ only involves colored forests with less than $k$ coalescent branches, from the original root, up to the final level. If $\|\mathbf{p}\| \geq (2k)$, then these colored forests contain at least one elementary tree with a white leaf. By definition of the operator $\overline{\mathbf{Q}}_n^{(\mathbf{p}+q)}$ we find that

$$\partial^k \overline{\mathbb{Q}}_{n,\mathbf{p}+q} \overline{\mathbf{Q}}_n^{(\mathbf{p}+q)}(F_{n+1}) = 0.$$

This yields the second part of the theorem. □

**5. Consequences.** We end this paper with a series of some direct consequences of the above theorem. Some of the already known properties of particle systems can be found in Section 1.

- The differential forest expansions presented in Theorem 4.12 allow us to deduce precise strong propagation of chaos estimates. For instance, for any $N \geq (q+7)$ [so that $(\lfloor (N-q)/2 \rfloor - 1) \geq 2$] we have

$$\sup_{N\geq q+7} N^2 \left\| \mathbb{P}_{n+1,q}^N - \eta_{n+1}^{\otimes q} - \frac{1}{N}\partial^1 \mathbb{P}_{n+1,q} \right\|_{\mathrm{TV}} < \infty$$

with a first-order partial derivative given by the formula

$$\partial^1 \mathbb{P}_{n+1,q}(F_{n+1}) = \partial^1 \mathbb{Q}_{n,q} Q_{n+1}^{\otimes q}(\overline{F}_{n+1}) + q \sum_{0\leq m\leq n} \partial^1 \overline{\mathbb{Q}}_{n,\mathbf{q}_m} \overline{\mathbf{Q}}_n^{\mathbf{q}_m}(F_{n+1})$$

with the sequence of integers $\mathbf{q}_m = (1_m(k) + q1_n(k))_{0\leq k\leq n}$. The first term in the right-hand side in the above displayed formula has been treated in Corollary 3.13, and we have that

$$\partial^1 \mathbb{Q}_{n,q} Q_{n+1}^{\otimes q}(\overline{F}_{n+1})$$
$$= \frac{q(q-1)}{2} \sum_{0\leq k\leq n} \Delta_{n,q}^{\mathbf{f}_{1,k}} Q_{n+1}^{\otimes q}(\overline{F}_{n+1})$$
$$= \frac{q(q-1)}{2} \sum_{0\leq k\leq n} \gamma_k(1) \int_{E_k^{q-1}} \gamma_k^{(q-1)}(d(x^2,\ldots,x^q))$$
$$\times Q_{k,n+1}^{\otimes q}(\overline{F}_{n+1})(x^2, x^2, x^3, \ldots, x^q).$$

Each of the terms $\partial^1 \overline{\mathbb{Q}}_{n,\mathbf{q}_m} \overline{\mathbf{Q}}_n^{\mathbf{q}_m}(F_{n+1})$ only involves the colored forests

$$\forall 0 \leq k \leq m \qquad \overline{\mathbf{f}}_{k,m} := U_k T_{k,m,n} U_{n+1}^{q-1}$$

associated with the trees $T_{k,m,n}$, and $U_l$, introduced on pages 816–817. After some elementary manipulations, we find that

$$\partial^1 \overline{\mathbb{Q}}_{n,\mathbf{q}_m} \overline{\mathbf{Q}}_n^{\mathbf{q}_m}(F_{n+1})$$

$$= q \sum_{0 \leq k \leq m} \Delta_{n,\mathbf{q}_m}^{\overline{\mathbf{f}}_{k,m}} \overline{\mathbf{Q}}_n^{\mathbf{q}_m}(F_{n+1})$$

$$= q \sum_{0 \leq k \leq m} \gamma_k(1) \int_{E_k^q} \gamma_k^{\otimes q}(d(x^1, \ldots, x^q)) Q_{k,m}(\overline{G}_m)(x^1)$$

$$\times Q_{k,n+1}^{\otimes q}(\overline{F}_{n+1})(x^1, \ldots, x^q).$$

- A Wick formula derives from Theorem 4.9. More precisely, for any $F \in \mathcal{B}_0^{\mathrm{sym}}(E_{n+1}^q)$, the partial derivatives $\partial^k \mathbb{P}_{n+1,q}(F)$ are null up to any order $k \leq \lfloor q/2 \rfloor$, and for any even integer $q$, we have

$$(5.1) \qquad \partial^{q/2} \mathbb{P}_{n+1,q}(F) = \gamma_{n+1}(1)^{-q} \partial^{q/2} \mathbb{Q}_{n,q} Q_{n+1}^{\otimes q}(F)$$

as soon as $N \geq 2(q+2)$ [so that $(\lfloor (N-q)/2 \rfloor - 1) \geq q/2$]. To prove this claim, we notice that for any $k < q/2$, and any $F \in \mathcal{B}_0^{\mathrm{sym}}(E_{n+1}^q)$, we have

$$\forall \mathbf{p} \in \bigcup_{0 \leq l < 2k} \mathbb{N}_l^{n+1}, \qquad k < q/2 \leq (\|\mathbf{p}\| + q)/2 \quad \text{and therefore}$$

$$\partial^k \overline{\mathbb{Q}}_{n,\mathbf{p}+q} \overline{\mathbf{Q}}_n^{(\mathbf{p}+q)}(F) = 0.$$

This yields that $\partial^k \mathbb{P}_{n+1,q}(F) = 0$, for any even integer $q$, and any $k < q/2$. In the case $k = q/2$, we have $k = q/2 = (\|\mathbf{p}\| + q)/2$ if, and only if, $\mathbf{p}$ coincides with the null sequence $\mathbf{0}$.

- We let $\|\mu\|_{\mathcal{B}_{0,1}^{\mathrm{sym}}} = \sup_{F \in \mathcal{B}_{0,1}^{\mathrm{sym}}} |\mu(F)|$ be the Zolotarev seminorm on $\mathcal{M}(E_{n+1}^q)$ associated with the collection of functions

$$\mathcal{B}_{0,1}^{\mathrm{sym}} := \{F \in \mathcal{B}_0^{\mathrm{sym}}(E_{n+1}^q) : \|F\| \leq 1\}.$$

For any even integer $q$ such that $(q+2) \leq N/2$, we have

$$(5.2) \qquad \sup_{N \geq q+7} N^{1+q/2} \left\| \mathbb{P}_{n+1,q}^N - \frac{1}{N^{q/2}} \frac{1}{\gamma_{n+1}(1)^q} \partial^{q/2} \mathbb{Q}_{n,q} Q_{n+1}^{\otimes q} \right\|_{\mathcal{B}_{0,1}^{\mathrm{sym}}} < \infty.$$

- Combining the Wick formula stated above with the Borel–Cantelli lemma, we obtain for all $q \geq 4$ the almost sure convergence result

$$(5.3) \qquad \lim_{N \to \infty} (\eta_n^N)^{\odot q}(F) = 0 \qquad \text{p.s.}$$

for any bounded symmetric function $F \in \mathcal{B}_0^{\mathrm{sym}}(E_n^q)$. This result is an extension of the law of large numbers for $U$-statistics obtained by Hoeffding [8] for independent and identically distributed random variables to interacting particle models. One can also look at [9] for a more modern reference.

- We mention that the same lines of arguments used in the proof of Theorem 4.12 show that the sequence of probability measures

$$\tilde{\mathbb{P}}^N_{n+1,q} : F \in \mathcal{B}_b(E^q_{n+1}) \mapsto \tilde{\mathbb{P}}^N_{n+1,q}(F) := \mathbb{E}((\eta^N_{n+1})^{\otimes q}(F))$$

is differentiable up to order $\lfloor (N-q)/2 \rfloor$, with $\partial^0 \tilde{\mathbb{P}}_{n+1,q} = \eta^{\otimes q}_{n+1}$, and the partial derivatives given for any $1 \leq k < \lfloor (N-q)/2 \rfloor$ by the following formula:

$$\partial^k \tilde{\mathbb{P}}_{n+1,q}(F) = \sum_{\mathbf{p} \in \bigcup_{0 \leq l < 2k} \mathbb{N}^{n+1}_l} \frac{((q-1)+|\mathbf{p}|)!}{(q-1)!\mathbf{p}!} \partial^k \overline{\mathbb{Q}}_{n+1,(\mathbf{p},q)}(\overline{G}^{\mathbf{p}}_n \otimes \overline{F}).$$

In the same way, for any $F \in \mathcal{B}^{\mathrm{sym}}_0(E^q_{n+1})$ and $N \geq 2(q+2)$, the partial derivatives $\partial^k \mathbb{Q}_{n+1,q}(F)$ are null up to any order $k < q/2$, and for any even integer $q$, we have that

$$\partial^{q/2} \tilde{\mathbb{P}}_{n+1,q}(F) = \gamma_{n+1}(1)^{-q} \partial^{q/2} \mathbb{Q}_{n+1,q}(F).$$

- Finally, the Wick formula stated above allows to deduce sharp $\mathbb{L}_q$-mean error bound. To see this claim, we simply observe that

$$F = (f - \eta_{n+1}(f))^{\otimes q}$$

with $f \in \mathcal{B}_b(E_{n+1}) \Longrightarrow \tilde{\mathbb{P}}^N_{n+1,q}(F) = \mathbb{E}([\eta^N_{n+1}(f) - \eta_{n+1}(f)]^q)$

and for any even integer $q$ such that $(q+2) \leq N/2$, we have the inequality

$$(5.4) \qquad \sup_{N \geq q+7} N^{1+q/2} \left\| \tilde{\mathbb{P}}^N_{n+1,q} - \frac{1}{N^{q/2}} \frac{1}{\gamma_{n+1}(1)^q} \partial^{q/2} \mathbb{Q}_{n+1,q} \right\|_{\mathcal{B}^{\mathrm{sym}}_{0,1}} < \infty.$$

## REFERENCES

[1] BROUDER, C. (2004). Trees, renormalization and differential equations. *BIT* **44** 425–438. MR2106008
[2] COLLINS, J. C. (1984). *Renormalization.* Cambridge Univ. Press. MR778558
[3] COMTET, L. (1970). *Analyse Combinatoire. Tomes I, II. Collection SUP: "Le Mathématicien"* **4–5**. Presses Universitaires de France, Paris. MR0262087
[4] CONNES, A. and KREIMER, D. (2000). Renormalization in quantum field theory and the Riemann–Hilbert problem. I. The Hopf algebra structure of graphs and the main theorem. *Comm. Math. Phys.* **210** 249–273. MR1748177
[5] DEL MORAL, P. (2004). *Feynman–Kac Formulae: Genealogical and interacting particle systems with applications.* Springer, New York. MR2044973
[6] DOUCET, A., DE FREITAS, N. and GORDON, N., eds. (2001). *Sequential Monte Carlo Methods in Practice.* Springer, New York. MR1847783
[7] GUIONNET, A. and SEGALA MAUREL, E. (2006). Combinatorial aspects of matrix models. *ALEA* **1** 241–279.
[8] HOEFFDING, W. (1961). The strong law of large numbers for $U$-statistics. Mimeo Report 302, Inst. Statist., Univ. of North Carolina.
[9] LEE, A. J. (1990). *$U$-Statistics: Theory and Practice. Statistics: Textbooks and Monographs* **110**. Marcel Dekker, New York. MR1075417

P. Del Moral
Centre INRIA Bordeaux Sud-Ouest
  & Institut de Mathématiques
  de Bordeaux
UMR 5251, CNRS
Université Bordeaux I
351, cours de la Libération
33405 Talence cedex
France
E-mail: delmoral@math.u-bordeaux1.fr

F. Patras
S. Rubenthaler
Laboratoire J.-A. Dieudonné
UMR 6621, CNRS
Parc Valrose
Université de Nice-Sophia
  Antipolis
06108 Nice cedex 02
France
E-mail: patras@math.unice.fr
        rubentha@unice.fr